\documentclass[conf]{new-aiaa}
\usepackage[utf8]{inputenc}

\usepackage{graphicx}
\usepackage{amsmath, amsfonts}
\usepackage[nameinlink]{cleveref}
\usepackage{siunitx}
\usepackage{longtable, tabularx}
\usepackage{algorithm}
\usepackage[noend]{algpseudocode}
\usepackage{subcaption}
\usepackage{pgfplots}
\pgfplotsset{compat=1.18}
\usepackage[dvipsnames]{xcolor}
\usepackage{soul}
\usepackage{tabto}

\definecolor{BYUblue}{RGB}{0, 61, 165}
\hypersetup{
    colorlinks=true,
    linkcolor=blue,
    citecolor=blue,
    urlcolor=blue,
    pdftitle={A Dynamic Subspace Approach for Low-rank Approximation of Large-scale Nonlinear Systems},
    pdfauthor={Jack DeChant, Rudy Geelen, Shane A McQuarrie, Johann Guilleminot}
}

\usepackage{silence}
\DeclareMathOperator*{\argmin}{arg\,min} 

\algrenewcommand\algorithmicrequire{\textbf{Input:}}
\algrenewcommand\algorithmicensure{\textbf{Output:}}
\algnewcommand{\LineComment}[1]{\vspace{.5em}\noindent\hspace{-1.1em}\textcolor{gray}{\texttt{\# #1}}}
\algrenewcomment[1]{\tabto{2.5in}\textcolor{gray}{\texttt{\# #1}}}

\crefname{equation}{}{}
\Crefname{equation}{}{}
\crefname{figure}{Fig.}{}
\Crefname{figure}{Fig.}{}

\title{A Dynamic Subspace Approach for Low-rank Approximation \\ of Large-scale Nonlinear Systems}

\author{Jack DeChant\footnote{Undergraduate Researcher} and Rudy Geelen\footnote{AIAA Member}}
\affil{Department of Aerospace Engineering, Texas A\&M University, College Station, TX, 77843}

\author{Shane~A.~McQuarrie}
\affil{Department of Mathematics, Brigham Young University, Provo, UT, 84602}

\author{Johann Guilleminot}
\affil{Department of Mechanical Engineering and Materials Science, Duke University, Durham, NC, 27708}

\begin{document}

\maketitle

\begin{abstract}
We present a dynamic subspace approach for efficiently approximating large-scale systems by learning time-continuous trajectories on the Grassmannian manifold. By parameterizing a low-dimensional basis as a geodesic path, the method allows for adaptive tracking of evolving physics. Our approach decouples the geometric drift of the subspace from the intrinsic state evolution. This avoids the typical rank inflation required by static low-dimensional approximation methods to maintain accuracy, effectively breaking the Kolmogorov barrier in transport-dominated phenomena. To ensure scalability for high-dimensional data, the optimization is performed in a reduced feature space, rendering the computational cost independent of the large original state dimension. Numerical results for a 1D transport equation and a large-scale turbulent airfoil wake demonstrate that this dynamic subspace approach achieves higher accuracy than static linear approximations at equivalent ranks, positioning it as a robust and scalable method for the low-rank modeling of complex, non-stationary dynamical systems.
\end{abstract}

\section{Nomenclature}

{\renewcommand\arraystretch{1.0}
\noindent\begin{longtable*}{@{}l @{\quad=\quad} l@{}}
$t$ & Time variable \\
$\tau$ & Normalized time variable \\
$N$ & Dimension of the full-order state vector \\
$K$ & Number of training snapshots \\
$r$ & Reduced state dimension \\
$n_f$ & Number of feature modes used for scalable optimization \\
$\mathbf{q}(t)$ & State vector, $\mathbf{q}(t) \in \mathbb{R}^N$ \\
$\bar{\mathbf{q}}$ & Reference state, $\bar{\mathbf{q}} \in \mathbb{R}^N$ \\
$\hat{\mathbf{q}}(t)$ & Generalized coefficient vector, $\hat{\mathbf{q}}(t) \in \mathbb{R}^r$ \\
$\mathbf{V}(\tau)$ & Time-varying basis matrix on the Grassmannian, $\mathbf{V}(\tau) \in \mathbb{R}^{N \times r}$ \\
$\boldsymbol{\Gamma}$ & Initial ``velocity'' matrix (tangent vector), $\boldsymbol{\Gamma} \in \mathbb{R}^{n \times r}$ \\
$\theta_1,\theta_2,\ldots,\theta_r$ & Principal angles, $\theta_i\in\mathbb{R}$ \\
$\boldsymbol{\Theta}$ & Diagonal matrix of principal angles, $\boldsymbol{\Theta} = \operatorname{diag}(\theta_1,\theta_2,\ldots,\theta_r)\in\mathbb{R}^{r\times r}$ \\
$\mathbf{V}_1,\mathbf{V}_2$ & Orthogonal rotation matrices in full state space, $\mathbf{V}_i\in\mathbb{R}^{N\times r}$ \\
$\boldsymbol{\Omega}$ & Joint full rotation matrix, $\boldsymbol{\Omega} = [~\mathbf{V}_1~~\mathbf{V}_2~]\in\mathbb{R}^{N\times 2r}$ \\
$\tilde{\mathbf{V}}$ & Matrix of feature modes (leading left singular vectors), $\tilde{\mathbf{V}} \in \mathbb{R}^{N \times n_f}$ \\
$\mathbf{T}_1, \mathbf{T}_2$ & Orthogonal rotation matrices in feature space, $\mathbf{T}_i \in \mathbb{R}^{n_f \times r}$ \\
$\mathbf{T}$ & Joint feature rotation matrix, $\mathbf{T} = [~\mathbf{T}_1~~\mathbf{T}_2~] \in \operatorname{St}(n_f, 2r)$ \\
$\mathbf{Q}$ & Snapshot matrix of observed states, $\mathbf{Q} \in \mathbb{R}^{N \times K}$ \\
$\mathbf{Y}$ & Feature space representation of the snapshots, $\mathbf{Y} \in \mathbb{R}^{n_f \times K}$ \\
$\mathbf{I}_a$ & $a \times a $ identity matrix \\
$\mathbf{C}$ & Cosines of the principal angles at $K$ time instances, $\mathbf{C} \in \mathbb{R}^{K\times r}$\\
$\mathbf{S}$ & Sines of the principal angles at $K$ time instances, $\mathbf{S} \in \mathbb{R}^{K\times r}$\\
$\odot$ & Hadamard (elementwise) product \\
$\|\cdot\|_F$ & Frobenius norm \\
$\operatorname{St}(m, n)$ & Stiefel manifold of $m \times n$ matrices with orthonormal columns \\
$\mathcal{M}$ & Product manifold $\operatorname{St}(n_f, 2r) \times \mathbb{R}^r$, the search space for the optimization problem \\
$\mathcal{G}(r, n)$ & Grassmannian manifold of $r$-dimensional subspaces in $\mathbb{R}^n$ \\
$\mathcal{T}_{\mathcal{V}}\mathcal{G}(r, n)$ & The tangent space of the Grassmannian manifold at the subspace $\mathcal{V}\in\mathcal{G}(r, n)$
\end{longtable*}}

\section{Introduction}

Effective low-dimensional approximation is a key component of many applications in data compression, machine learning, and model reduction~\cite{benner2015pmorsurvey}. Common low-dimensional approximation methodologies,
ranging from classical Proper Orthogonal Decomposition (POD)~\cite{sirovich1987turbulence, berkooz1993proper} to modern neural autoencoders~\cite{lee2020model}, typically seek a \emph{stationary} approximation space or latent coordinate system. However, such static representations are often inadequate for capturing transient, transport-dominated phenomena in which physical modes or coherent structures drift over time. To address this lack of temporal adaptability and the resulting geometric misalignment, we propose a mathematical framework that replaces fixed projection paradigms with \emph{dynamic} subspaces.
By utilizing time-continuous trajectories on the Grassmannian manifold, we reframe dimensionality reduction as a differential-geometric task, which allows the approximation space to evolve alongside the high-dimensional state being described. This approach disambiguates the geometric drift of the subspace from the intrinsic evolution of the low-dimensional coordinates, ensuring the approximation remains synchronized with the instantaneous physics and providing a robust, scalable framework for low-dimensional modeling of complex systems that exceed the capabilities of stationary approximation spaces.

Static low-dimensional approximation methods involve a fundamental tradeoff between accuracy and interpretability. Linear or affine approximations such as POD, balanced truncation~\cite{gugercin2004balancedtruncation,gosea2022datadrivenbt}, and the reduced basis method~\cite{hesthaven2022reduced} build a linear subspace to capture the dominant behavior of the trajectory being approximated. Linear methods are highly interpretable, tend to be computationally straightforward, and are well suited for applications such as projection-based model reduction; on the other hand, many nonlinear systems cannot be well represented in a linear subspace unless the dimensionality of that subspace is large, diminishing any efficiency gains provided by the approximation. This limitation is known as the Kolmogorov barrier in the literature~\cite{peherstorfer2022breaking}. Nonlinear approximation methods can overcome this fundamental expressivity limitation, and thereby mitigate the Kolmogorov barrier, at the expense of some interpretability and/or computational efficiency. See, for instance,~\cite{barnett2022quadratic,geelen2023operator,geelen2024learning,schwerdtner2024greedyquad} for quadratic approximations,~\cite{scholkopf1997kernel,diaz2025kernelmanifold} for kernel-based methods,~\cite{lee2020model,fries2022lasdi} for model reduction applications with neural autoencoders, and~\cite{peherstorfer2022breaking} for a note on the Kolmogorov barrier as it relates to nonlinear model reduction. These and related machine learning techniques offer improved data compression but rely fundamentally on fixed embeddings that do not adapt to the state as time evolves. In this work, we consider time-continuous subspace approximation methods, which represent the evolving state as linear combinations of basis functions that evolve with time. The motivation behind this structure is that nonlinear processes can often be approximated well in linear subspaces \emph{locally in time}, even when approximations in a single, global reduced space are inefficient~\cite{hesthaven2026nonlinearmodelreductiontransportdominated}. Time-continuous linear approximation therefore addresses the limitations of static linear approximation methods while remaining highly interpretable.

There are several existing strategies across many application domains for introducing time dependence into low-dimensional approximations.
Early foundations for this idea can be found in the work on finite element methods with moving nodes~\cite{miller1981moving}, which represents one of the first Galerkin frameworks where both the coefficients and the basis functions are permitted to vary simultaneously in time. In the signal processing community, there exists a rich history of dynamic subspace \emph{tracking} methods designed to follow non-stationary data streams in real time~\cite{yang1995projection}. Other perspectives are provided by dynamic low-rank approximations which parametrize the solution as a low-rank matrix and evolve its factors over time~\cite{koch2007dynamical}. These approaches rely on the Dirac--Frenkel variational principle to project full dynamics onto the tangent space of a fixed-rank matrix manifold. This concept is closely related to the dynamic orthogonal decomposition developed for reduced-order modeling of stochastic differential systems, where the low-rank structure is extracted directly from the governing operators to capture complex transient or stochastic behavior~\cite{sapsis2009dynamically, babaee2016minimization, patil2020real, ramezanian2021fly}.
In the model reduction community, one prominent approach focuses on mimicking transport dynamics to evolve the reduced space by utilizing techniques such as optimal transport~\cite{iollo2014advection}, approximated Lax pairs~\cite{gerbeau2014approximated}, or characteristic-based transport~\cite{rim2023manifold}. Alternatively, time-dependent modes can be evolved through instantaneous residual minimization~\cite{black2020projection} or by periodically reverting to the full model to generate new snapshots~\cite{rapun2010reduced}. Another perspective is to perform efficient low-rank updates to the reduced-order model and the approximation space in a way that integrates real-time information, such as full-order model residuals or incomplete state fields, without a complete reconstruction of the reduced model. This is exemplified by the work of Peherstorfer and Willcox for linear and nonlinear problems~\cite{peherstorfer2015dynamic, peherstorfer2015online}. For a comprehensive overview of these and other model reduction methodologies, we recommend~\cite{hesthaven2022reduced}. A primary challenge for these methods lies with the classical offline-online decomposition of computational tasks, as adapting basis functions at runtime often incurs costs that scale with the high-dimensional state-space and diminishes speedups. To mitigate this, one alternative is to precompute time-dependent subspaces during the offline phase from training trajectories, allowing for online deployment without further adaptation or enrichment~\cite{billaud2017dynamical}.

Building upon the blueprint from~\cite{billaud2017dynamical} of utilizing \emph{a priori} computed dynamic low-rank approximations, we propose a robust data-driven framework that models evolving subspaces as geodesic trajectories on the Grassmannian manifold. By representing the shortest paths between fixed-dimensional subspaces, these geodesics allow the approximation to synchronize seamlessly with the system's instantaneous physics, maintaining a low-rank and spatially compact footprint that effectively avoids the rank inflation typically needed with stationary approximation methods to maintain accuracy. Crucially, this methodology remains entirely independent of the governing equations and requires no explicit parametrization of the underlying physics, making the approach straightforward to integrate into existing production-level simulation codes. By reframing dimensionality reduction as a formal differential-geometric task, we shift the computational focus away from the manipulation of high-dimensional basis vectors and toward the optimization of temporal coefficients and principal angles. This shift offers a mathematically rigorous, scalable alternative to traditional manifold learning as an enabler for high-fidelity digital twins in the exascale computing era.

This article is organized as follows. \Cref{sec:dynamic_space_estimation} introduces the theoretical foundation of dynamic subspace estimation, revisiting classical dimensionality reduction and detailing the differential geometry of the Grassmannian used to model geodesic trajectories. \Cref{sec:scalable_optimization} describes the scalable optimization framework, including the elimination of high-dimensional costs and the vectorized implementation within an intrinsic manifold-constrained environment. Numerical results are then presented in \Cref{sec:results}, where the framework is evaluated on a one-dimensional transport equation and a large-scale turbulent airfoil wake large eddy simulation (LES). Finally, concluding remarks and potential avenues for future work are provided in \Cref{sec:conlusions_and_future_work}.

\section{Dynamic Subspace Estimation}
\label{sec:dynamic_space_estimation}

\subsection{Time-continuous Low-rank Approximation}

Consider a high-dimensional time-dependent vector $\mathbf{q}(t) \in \mathbb{R}^{N}$ representing the state variable(s) of a physical system (e.g., temperature and/or pressure) with $N$ degrees of freedom.
Static linear approximation methods, such as Proper Orthogonal Decomposition (POD), find a fixed basis matrix $\mathbf{V}$ in the Stiefel manifold $\operatorname{St}(N, r)$, the set of all $N \times r$ matrices with orthonormal columns, to induce a state approximation
\begin{equation}
    \mathbf{q}(t) \approx \bar{\mathbf{q}} + \mathbf{V} \hat{\mathbf{q}}(t),
    \label{eq:pod}
\end{equation}
where $\bar{\mathbf{q}} \in \mathbb{R}^{N}$ is a fixed reference state (usually a time-averaged state or an initial condition) and $\hat{\mathbf{q}}(t) \in \mathbb{R}^r$ is a generalized coordinate vector.
When $r \ll N$, \cref{eq:pod} is a low-dimensional approximation of the state $\mathbf{q}(t)$ that depends entirely on the time-dependent coordinates $\hat{\mathbf{q}}(t)$.
Given $\bar{\mathbf{q}}$, the goal is to choose $\mathbf{V}$ so that $\mathbf{q}(t) - \bar{\mathbf{q}}$ lies in or near span of the columns of $\mathbf{V}$, a linear subspace in $\mathbb{R}^{N}$, for all $t$ in the time domain of interest. If this is the case, then $\hat{\mathbf{q}}(t)$ can be chosen so that the approximation error $\|\mathbf{q}(t) - \bar{\mathbf{q}} - \mathbf{V}\hat{\mathbf{q}}(t)\|_2$ is small for all $t$.
However, many time-varying nonlinear phenomena cannot be well represented within a fixed low-dimensional linear subspace, meaning a moderately large approximation rank $r$ is required to drive down the approximation error.

In this work, we increase the flexibility of the static low-dimensional approximation \cref{eq:pod} by allowing the basis matrix $\mathbf{V}$ to vary in time.
We therefore consider a \emph{dynamic} state approximation
\begin{equation}
    \mathbf{q}(t)
    \approx \bar{\mathbf{q}} + \mathbf{V}(t) \hat{\mathbf{q}}(t),
    \label{eq:dynamic}
\end{equation}
with the requirement that $\mathbf{V}(t)\in\operatorname{St}(N, r)$ for all $t$.
In this formulation, the columns of $\mathbf{V}(t)$ are dynamic basis vectors that adapt to the evolving state, allowing the reduced coordinates $\hat{\mathbf{q}}(t)$ to effectively characterize latent state dynamics within a coordinate frame that remains synchronized with the physical state.
Our objective is to parameterize and construct $\mathbf{V}(t)$ so that the approximation error $\|\mathbf{q}(t) - \bar{\mathbf{q}} - \mathbf{V}(t)\hat{\mathbf{q}}(t)\|_2$ can be controlled for all $t$ without inflating the approximation rank $r$.
For modeling physical phenomena it is also desirable that $\mathbf{V}(t)$ evolve continuously in time, for which we appeal to differential geometry.

\subsection{Differential Geometry of the Grassmannian Manifold}

\begin{figure}[tbp]
\centering
    \includegraphics[width=0.45\linewidth]{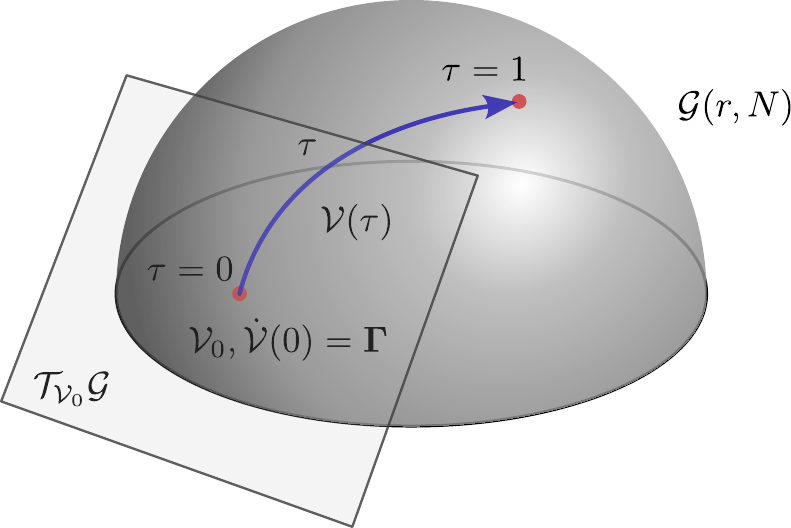}
    \caption{A time-continuous trajectory $\mathcal{V}(\tau)$ on the Grassmann manifold $\mathcal{G}(r,N)$ as defined by a single geodesic, represented by the blue curve, defined by the end points (red). The image illustrates how the matrix representatives $\mathbf{V}_0$ and $\mathbf{\Gamma}$ (constrained by $\mathbf{V}_0^\top \mathbf{\Gamma} = \mathbf{0}$) define the abstract path on the manifold.}
    \label{fig:geodesic_model}
\end{figure}

The Grassmannian $\mathcal{G}(r,N)$ is the smooth manifold consisting of all $r$-dimensional linear subspaces of $\mathbb{R}^N$.
An element $\mathcal{V}\in\mathcal{G}(r,N)$ of the Grassmannian can be efficiently represented numerically by an orthonormal matrix $\mathbf{V} = [~\mathbf{v}_1~~\cdots~~\mathbf{v}_r~] \in \operatorname{St}(N, r)$ such that $\operatorname{span}(\{\mathbf{v}_1,\ldots,\mathbf{v}_r\}) = \mathcal{V}$, that is, the columns of $\mathbf{V}$ form an orthonormal basis for $\mathcal{V}$.
To model the smooth evolution of a dynamic reduced-order basis $\mathbf{V}(t)$ in \cref{eq:dynamic}, we consider a \emph{geodesic} $\mathcal{V}(\tau)$, the shortest path between two points on the Grassmannian, with normalized time variable $\tau\in[0,1]$. As established in~\cite{edelman1998geometry, absil2004riemannian}, a geodesic $\mathcal{V}(\tau)$ is uniquely determined by its initial position $\mathcal{V}(0) = \mathcal{V}_0$ and an initial velocity $\dot{\mathcal{V}}(0)$, residing in the tangent space $\mathcal{T}_{\mathcal{V}_0} \mathcal{G}(r,N)$, which dictates the rotation of the subspace as it progresses along the geodesic.
For computational treatment, these abstract quantities are represented by the pair $(\mathbf{V}_0, \mathbf{\Gamma})$, where $\mathbf{V}_0 \in \operatorname{St}(N, r)$ corresponds to the initial position $\mathcal{V}_0$ and $\mathbf{\Gamma} \in \mathbb{R}^{N\times r}$ is the initial velocity matrix.
The set of matrices corresponding to the tangent space $\mathcal{T}_{\mathcal{V}_0} \mathcal{G}(r,N)$ is given by $\{\mathbf{Y}\in\mathbb{R}^{N\times r}:\mathbf{V}_0^\top\mathbf{Y} = \mathbf{0}\}$, hence $\mathbf{V}_0^\top \mathbf{\Gamma} = \mathbf{0}$. This requirement ensures that $\mathbf{\Gamma}$ represents a change in the subspace rather than a redundant rotation within the basis.

A numerical representation for the continuous evolution of the dynamic subspace $\mathcal{V}(\tau)$, that is, a matrix trajectory $\mathbf{V}(\tau)$ mapping $\tau\in[0,1]$ to an orthonormal matrix in $\operatorname{St}(N, r)$, satisfies the following second-order system of ordinary differential equations that arises from the calculus of variations for minimizing path length on $\mathcal{G}(r,N)$:
\begin{equation}
    \ddot{\mathbf{V}}(\tau)
    + \mathbf{V}(\tau) \left( \dot{\mathbf{V}}(\tau)^\top \dot{\mathbf{V}}(\tau) \right)
    = \mathbf{0},
    \qquad
    \mathbf{V}(0) = \mathbf{V}_0,
    \qquad
    \dot{\mathbf{V}}(0) = \mathbf{\Gamma}.
    \label{eq:ode}
\end{equation}
\Cref{fig:geodesic_model} provides a schematic overview of these variables, visualizing how the matrix representations and tangent space elements define the abstract path on the Grassmannian. The analytical solution to \cref{eq:ode} is obtained via the so-called exponential map, which projects the tangent vector $\mathbf{\Gamma}$ along the geodesic.
Letting $\boldsymbol{\Gamma} = \boldsymbol{\Psi\Sigma\Phi}^\top$ be a thin singular value decomposition (SVD) of the initial velocity matrix, $\mathbf{V}(t)$ has the closed-form trigonometric representation
\begin{equation}
    \mathbf{V}(\tau)
    = \mathbf{V}_0 \boldsymbol{\Phi} \cos(\tau \boldsymbol{\Sigma})
    + \boldsymbol{\Psi} \sin(\tau \boldsymbol{\Sigma}),
\label{eq:geodesic_model_equation}
\end{equation}
where $\cos$ and $\sin$ denote the matrix cosine and matrix sine (akin to the matrix exponential), respectively; because $\tau\boldsymbol{\Sigma}$ is diagonal, $\cos(\tau\boldsymbol{\Sigma}) = \operatorname{diag}(\cos(\tau\sigma_1),\cos(\tau\sigma_2),\ldots,\cos(\tau\sigma_r))$ in which $\Sigma = \operatorname{diag}(\sigma_1, \sigma_2, \cdots, \sigma_r)$, and likewise for $\sin(\tau\boldsymbol{\Sigma})$.
Geometrically, the diagonal entries of $\boldsymbol{\Sigma}$ represent the principal angles of rotation per unit time.

\subsection{Geodesic Model for Subspace Tracking}

From the above differential-geometric foundation, we adopt the structural form \cref{eq:geodesic_model_equation} in a way that is amenable to optimization and data-driven learning.
The main idea is to treat learning the geodesic as a geometric regression problem where the initial basis, the direction of evolution, and the associated principal angles are inferred from data produced by a physical model.
This closely resembles the approach taken in~\cite{blocker2023dynamicsubspaceestimationgrassmannian}; in \Cref{sec:scalable_optimization}, we develop a novel implementation strategy that scales well to high-dimensional data.

Consider a matrix representation for a continuous Grassmannian geodesic $\mathcal{V}(\tau)$ parameterized as
\begin{equation}
    \label{eq:cossin}
    \mathbf{V}(\tau)
    = \mathbf{V}_1 \cos(\tau \boldsymbol{\Theta})
    + \mathbf{V}_2 \sin(\tau \boldsymbol{\Theta}),
\end{equation}
where the columns of $\mathbf{V}_1 \in \mathbb{R}^{N \times r}$ form an orthonormal basis for a subspace $\mathcal{V}_1 \in \mathcal{G}(r,N)$ and $\mathbf{V}_2 \in \mathbb{R}^{N \times r}$ is an orthonormal matrix whose column span lies in the tangent space $\mathcal{T}_{\mathcal{V}_1}\mathcal{G}(r,N)$.
The diagonal matrix $\boldsymbol{\Theta} = \operatorname{diag}(\theta_1,\theta_2,\ldots,\theta_r) \in \mathbb{R}^{r \times r}$ contains the principal angles representing the fundamental rotations required to align the $r$-dimensional subspaces $\mathcal{V}(0)$ and $\mathcal{V}(1)$ as the normalized time $\tau$ increases from $\tau = 0$ to $\tau = 1$.
We then define the matrices
\begin{equation}
    \label{eq:omegaW}
    \boldsymbol{\Omega}
    := [~\mathbf{V}_1~~\mathbf{V}_2~] \in \mathbb{R}^{N \times 2r},
    \qquad
    \mathbf{W}(\boldsymbol{\Theta}, \tau)
    := \left[\begin{array}{c}
        \cos(\tau \boldsymbol{\Theta}) \\ \sin(\tau \boldsymbol{\Theta})  
    \end{array}\right]
    \in \mathbb{R}^{2r \times r},
\end{equation}
to obtain the following compact expression for \cref{eq:cossin}:
\begin{equation}
    \label{eq:fixed-dynamic-split}
    \mathbf{V}(\tau) = \boldsymbol{\Omega} \mathbf{W}(\boldsymbol{\Theta},\tau)
\end{equation}
Written this way, $\mathbf{V}(\tau)$ is completely determined by $\boldsymbol{\Omega}$ and the principal angles $\theta_1,\theta_2,\ldots,\theta_r$.
Substituting \cref{eq:fixed-dynamic-split} into \cref{eq:dynamic} and using the normalized time $\tau = \tau(t)$ yields a state approximation of the form
\begin{equation}
\begin{aligned}
    \mathbf{q}(t)
    \approx \bar{\mathbf{q}} + \mathbf{V}(\tau) \hat{\mathbf{q}}(t)
    &= \bar{\mathbf{q}} + \boldsymbol{\Omega} \mathbf{W}(\boldsymbol{\Theta},\tau) \hat{\mathbf{q}}(t).
    \label{eq:model_form}
\end{aligned}
\end{equation}
Our objective now is to identify $\boldsymbol{\Omega}$ and $\boldsymbol{\Theta}$ so that \cref{eq:model_form} is an accurate approximation for all $t$.

Given an initial observation time $t_1$ and a final observation time $t_K$, let $t_1 < t_2 < \cdots < t_K$ be discrete times for which data for the state $\mathbf{q}(t)$ are available for computation.
The normalized time is then $\tau(t) = (t - t_0) / (t_f - t_0)$, and we define $\tau_j = \tau(t_j)$ and denote the state observations, also called snapshots, with $\mathbf{q}(t_j)\in\mathbb{R}^{N}$ for $j=1,2,\ldots,N$.
For a given snapshot $\mathbf{q}(t_j)$, the reduced coordinates $\hat{\mathbf{q}}(t_j) = \mathbf{V}(\tau_j)^\top(\mathbf{q}(t_j) - \bar{\mathbf{q}})$ give the best possible approximation of $\mathbf{q}(t_j)$ through \cref{eq:model_form} in the $2$-norm sense.
In this case, the approximation error
\begin{equation}
    \begin{aligned}
    \|\mathbf{q}(t_j) - \bar{\mathbf{q}} - \mathbf{V}(\tau_j)\hat{\mathbf{q}}(t_j)\|_2
    &=
    \|\mathbf{q}(t_j) - \bar{\mathbf{q}} - \mathbf{V}(\tau_j)\mathbf{V}(\tau_j)^\top(\mathbf{q}(t_j) - \bar{\mathbf{q}})\|_2
    \\
    &=
    \|(\mathbf{I} - \mathbf{V}(\tau_j)\mathbf{V}(\tau_j)^\top)(\mathbf{q}(t_j) - \bar{\mathbf{q}})\|_2
    \end{aligned}
\end{equation}
is also called the projection error.
Our strategy is to build the dynamic subspace by minimizing the mean squared projection error, meaning we solve the optimization problem
\begin{equation}
    \begin{aligned}
    \min_{\boldsymbol{\Omega}, \boldsymbol{\Theta}} \frac{1}{K}&\sum_{j=1}^K \left\|
        \mathbf{q}(t_j)
        - \bar{\mathbf{q}}
        - \mathbf{V}(\tau_j) \mathbf{V}(\tau_j)^\top  (\mathbf{q}(t_j)
        - \bar{\mathbf{q}})
    \right\|_2^2
    \end{aligned}    
    \quad \text{subject to} \quad
    \boldsymbol{\Omega}^\top \boldsymbol{\Omega} = \mathbf{I}_{2r}.
\label{eq:optimization}
\end{equation}
The orthogonality constraint on $\boldsymbol{\Omega}$ ensures that $\mathbf{V}_1$ and $\mathbf{V}_2$ are each orthonormal and mutually orthogonal.
Because $\boldsymbol{\Theta}$ is diagonal, the total number of unknowns is $r(2N + 1)$.

It should be noted that the concept of normalized time $\tau$ is employed strictly for the purposes of parametrizing the geodesic. However, the methodology remains effective even when the available data lack strict temporal coherence, making the framework equally suited for broader parametric problems where the subspace evolution is mapped to system parameters. Here, we focus on temporal problems to demonstrate the ability to track time-evolving physics.

\section{Scalable Optimization Framework}
\label{sec:scalable_optimization}

\subsection{Eliminating High Dimensionality} 

The primary challenge in solving optimization problem~\cref{eq:optimization} stems from the state dimension $N$.
For large-scale systems, where (for example) $N$ may be on the order of millions, solving \cref{eq:optimization} directly is computationally intractable since it involves searching for an optimal ${N\times 2r}$ matrix $\boldsymbol{\Omega}$.
To address this challenge, we make the modeling assumption that the entire geodesic path resides within, or can be closely approximated by, the span of a fixed set of $n_f \ll N$ feature vectors, organized into an orthonormal matrix $\tilde{\mathbf{V}} \in \mathbb{R}^{N \times n_f}$.
This assumption is generally valid in practice if the feature matrix is constructed from a sufficiently rich ensemble of state snapshots by an information-maximizing procedure.
Specifically, we construct $\tilde{\mathbf{V}}$ using (centered) POD: organizing the snapshots into the data matrix $\mathbf{Q} = [~\mathbf{q}(t_1)~~\cdots~~\mathbf{q}(t_K)~]\in\mathbb{R}^{N\times K}$, we define $\tilde{\mathbf{V}}$ as the matrix whose columns are the first $n_f$ left singular vectors of the centered data matrix $\mathbf{Q} - \bar{\mathbf{q}}\mathbf{1}^\top = [~(\mathbf{q}(t_1) - \bar{\mathbf{q}})~~\cdots~~(\mathbf{q}(t_K) - \bar{\mathbf{q}})~]$.
The dimension $n_f$ is selected so that the resulting mean squared POD projection error
\begin{equation}
    \frac{1}{K}\sum_{j=1}^{K}\left\|
        \big(\mathbf{I} - \tilde{\mathbf{V}}\tilde{\mathbf{V}}^\top\big)
        (\mathbf{q}(t_j) - \bar{\mathbf{q}})
    \right\|_2^2
    \label{eq:pod_loss}
\end{equation}
is negligible, which is typically accomplished for $r \ll n_f \ll N$.
Thus, the feature space is large enough to accurately represent $\mathbf{q}(t)$ for all $t$, but small enough to enable (as we will see shortly) a significant computational speedup.

To constrain the search space to this precomputed subspace, we set $\mathbf{V}_1 = \tilde{\mathbf{V}}\mathbf{T}_1$ and $\mathbf{V}_2 = \tilde{\mathbf{V}}\mathbf{T}_2$
where $\mathbf{T}_1, \mathbf{T}_2 \in \mathbb{R}^{n_f \times r}$. The state approximation \cref{eq:model_form} then becomes
\begin{equation} 
\begin{aligned}
    \mathbf{q}(t)
    &\approx \bar{\mathbf{q}} + (\tilde{\mathbf{V}} \mathbf{T}_1 \cos(\tau \boldsymbol{\Theta})
    + \tilde{\mathbf{V}} \mathbf{T}_2 \sin(\tau \boldsymbol{\Theta})) \hat{\mathbf{q}}(t)
    \\
    &= \bar{\mathbf{q}} + \tilde{\mathbf{V}} \mathbf{T} \mathbf{W}(\boldsymbol{\Theta}, \tau) \hat{\mathbf{q}}(t)
\label{eq:approximation}
\end{aligned}
\end{equation}
where 
\begin{equation}
    \mathbf{T} = [~\mathbf{T}_1~~\mathbf{T}_2~] \in \mathbb{R}^{n_f \times 2r}.
    \label{eq:T}
\end{equation}
Because $\tilde{\mathbf{V}}$ is orthonormal, the matrices $\mathbf{V}_1$ and $\mathbf{V}_2$ are mutually orthonormal whenever $\mathbf{T}_1$ and $\mathbf{T}_2$ are mutually orthonormal.
Specifically, for $i,j\in\{1,2\}$,
\begin{align}
    \mathbf{V}_i^\top \mathbf{V}_j
    = (\tilde{\mathbf{V}} \mathbf{T}_i)^\top \tilde{\mathbf{V}} \mathbf{T}_j
    = \mathbf{T}_i^\top \tilde{\mathbf{V}}^\top \tilde{\mathbf{V}} \mathbf{T}_j
    = \mathbf{T}_i^\top \mathbf{T}_j
    = \begin{cases}
        \mathbf{I}_r &\text{if}~i=j, \\
        \mathbf{0} &\text{if}~i\neq j,
    \end{cases}
\end{align}
and hence $\boldsymbol{\Omega}^\top\boldsymbol{\Omega} = \mathbf{I}_{2r}$ if and only if $\mathbf{T}^\top\mathbf{T} = \mathbf{I}_{2r}$.
Thus, introducing $\tilde{\mathbf{V}}$ and enforcing orthonormality in $\mathbf{T}_1$ and $\mathbf{T}_2$ preserves the key geometric conditions of the Grassmannian trajectory.

The above reformulation greatly reduces the computational expense of the training task.
Instead of solving for the large orthonormal matrix $\boldsymbol{\Omega}\in\mathbb{R}^{N\times 2r}$, we optimize for orthonormal $\mathbf{T}\in\mathbb{R}^{n_f\times 2r}$ and the principal angles $\theta_1,\theta_2,\ldots,\theta_r$ by solving
\begin{equation}
    \min_{\mathbf{T}, \boldsymbol{\Theta} } \frac{1}{K}\sum_{j=1}^K \left\|
        \mathbf{q}(t_j)
        - \bar{\mathbf{q}}
        - \mathbf{V}(\tau_j) \mathbf{V}(\tau_j)^\top  (\mathbf{q}(t_j)
        - \bar{\mathbf{q}})
    \right\|_2^2
    \quad
    \text{subject to} \quad \mathbf{T}^\top \mathbf{T} = \mathbf{I}_{2r},
    \label{eq:reformulated0}
\end{equation}
or, more explicitly,
\begin{equation}
    \min_{\mathbf{T}, \boldsymbol{\Theta} } \frac{1}{K}\sum_{j=1}^K \left\|
        \big(\mathbf{I} - \tilde{\mathbf{V}}\mathbf{T}\mathbf{W}(\boldsymbol{\Theta},\tau_j)\mathbf{W}(\boldsymbol{\Theta},\tau_j)^\top\mathbf{T}^\top\tilde{\mathbf{V}}^\top\big)
        (\mathbf{q}(t_j) - \bar{\mathbf{q}})
    \right\|_2^2
    \quad
    \text{subject to} \quad \mathbf{T}^\top \mathbf{T} = \mathbf{I}_{2r}.
    \label{eq:reformulated}
\end{equation}
By shifting the focus to a coordinate-based representation within the feature space, the number of unknowns is reduced to $r(2n_f + 1)$ down from $r(2N + 1)$. This effectively decouples the optimization cost from the state dimension $N$. When the number of features $n_f$ is modest, this approach ensures that the numerical optimization remains fast and robust without sacrificing the ability to track complex, non-stationary dynamics.

\subsection{Feature Space Formulation and Vectorized Implementation}
\label{subsec:vectorization}

The numerical solution to the reformulated optimization problem \cref{eq:reformulated} is achieved through a manifold-constrained framework. Given that the subspace evolves continuously in time, a na\"ive implementation would require iterating through each discrete time step to update the basis $\mathbf{V}(\tau_j)$ and compute the corresponding projections. However, in the following, we show that all projections may be computed efficiently in a single pass. This effectively eliminates the computational bottleneck and overhead typically associated with looping over time steps, allowing the training phase to remain efficient even when using a time-continuous basis.

Because the feature matrix $\tilde{\mathbf{V}}$ in \cref{eq:approximation} is orthonormal, the $\ell_2$ norm on the high-dimensional state space equals the $\ell_2$ norm in the reduced feature space, that is, $\|\tilde{\mathbf{V}}\mathbf{x}\|_2 = \|\mathbf{x}\|_2$ for all $\mathbf{x}\in\mathbb{R}^{n_f}$.
This establishes a computationally efficient objective by allowing us to focus on $n_f$-dimensional coordinates.
To see this, we first project the data snapshots into the feature space as
\begin{equation}
    \mathbf{Y}
    = \tilde{\mathbf{V}}^\top(\mathbf{Q} - \bar{\mathbf{q}}\mathbf{1}^\top)
    = [~\mathbf{y}(t_1)~~\cdots~~\mathbf{y}(t_K)~]  \in \mathbb{R}^{n_f \times K},
    \label{eq:feature_space_representation}
\end{equation}
where
\begin{equation}
    \mathbf{y}_j = \mathbf{y}(t_j)
    = \tilde{\mathbf{V}}^\top (\mathbf{q}(t_j) - \bar{\mathbf{q}}) \in \mathbb{R}^{n_f}.
    \label{eq:feature_space_representation2}
\end{equation}
We also define
\begin{equation}
    \mathbf{R}(\mathbf{T, \boldsymbol{\Theta}},\mathbf{Y})
    = \left[\begin{array}{cc} 
        (\mathbf{C}^\top \odot \mathbf{C}^\top) \odot (\mathbf{T}_1^\top \mathbf{Y}) +
        (\mathbf{C}^\top \odot \mathbf{S}^\top) \odot (\mathbf{T}_2^\top \mathbf{Y})
        \\ 
        (\mathbf{S}^\top \odot \mathbf{C}^\top) \odot (\mathbf{T}_1^\top \mathbf{Y}) +
        (\mathbf{S}^\top \odot \mathbf{S}^\top) \odot (\mathbf{T}_2^\top \mathbf{Y})
    \end{array}\right]
    \in \mathbb{R}^{2r \times K},
\label{eq:R}
\end{equation}
where $\odot$ denotes the Hadamard (elementwise) product and the rotation components $\mathbf{C} = \mathbf{C}(\boldsymbol{\Theta})\in\mathbb{R}^{K\times r}$ and $\mathbf{S} = \mathbf{S}(\boldsymbol{\Theta})\in\mathbb{R}^{K\times r}$ have elements $(\mathbf{C})_{ji} = \cos(\tau_j\theta_i)$ and $(\mathbf{S})_{ji} = \sin(\tau_j \theta_i)$, respectively.
Then the approximation of the snapshot data in the feature space may be written, in matrix form, as
\begin{equation}
    \mathbf{Y} \approx \mathbf{T} \mathbf{R}(\mathbf{T, \boldsymbol{\Theta}},\mathbf{Y}).
    \label{eq:approx_matrix_form}
\end{equation}
A detailed derivation of feature space formulation and vectorization to obtain \cref{eq:R}--\cref{eq:approx_matrix_form} is deferred to the Appendix.

With the above definitions in hand, we can succinctly define the optimization problem in the feature space as
\begin{equation}
    \min_{(\mathbf{T}, \boldsymbol{\Theta}) \in \mathcal{M}}
    \frac{\|\mathbf{Y} - \mathbf{T} \mathbf{R}(\mathbf{T}, \boldsymbol{\Theta}, \mathbf{Y})\|_F}{\| \mathbf{Y}\|_F}
    \label{eq:riemannian_optimization}.
\end{equation}
To solve this manifold-constrained optimization problem, we utilize the \texttt{Pymanopt} package~\cite{townsend2016pymanopt}, which treats the search space as the product manifold $\mathcal{M} = \operatorname{St}(n_f, 2r) \times \mathbb{R}^r$, where $\operatorname{St}(n_f, 2r)$ is the Stiefel manifold of $n_f \times 2r$ orthonormal matrices and $\mathbb{R}^r$ is the Euclidean space for the principal angles $\boldsymbol{\Theta}$.
Rather than promoting orthogonality through traditional penalty functions or Lagrange multipliers, \texttt{Pymanopt} employs a Riemannian conjugate gradient solver to perform an optimization search that is strictly confined to $\mathcal{M}$.
Within this framework, the Riemannian gradient is calculated by projecting the Euclidean gradient onto the tangent space of $\mathcal{M}$, while a retraction mapping is used to wrap updated tangent vectors back onto $\mathcal{M}$ to ensure the basis remains strictly orthonormal at every iteration. 
For the conjugate gradient algorithm, vector transport is employed to shift previous search directions into the current tangent space, allowing the solver to combine directional information across the curved geometry. Our implementation leverages the \texttt{JAX} library~\cite{jax2018github} as the computational backend for high-performance automatic differentiation and just-in-time compilation to map exact derivatives directly to the tangent space.

\subsection{Motivational Example}

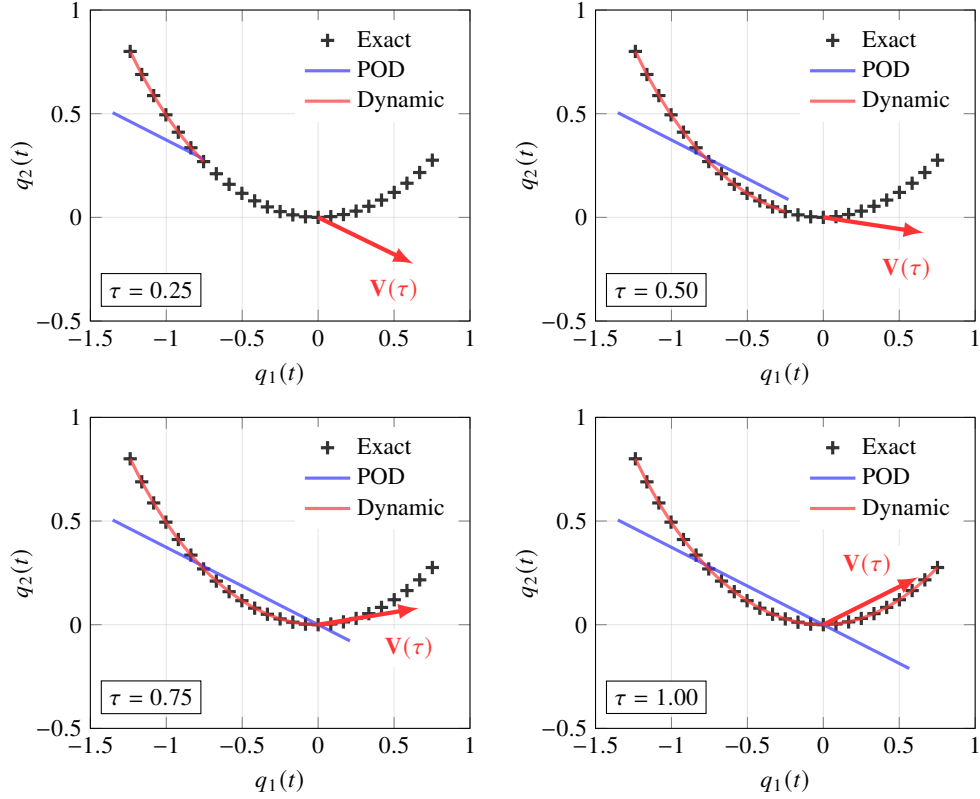
\begin{figure}[!tbp]
    \centering \small
    \begin{subfigure}[b]{0.4\linewidth}
    \centering
    \begin{tikzpicture}
    \begin{axis}[
        width=\linewidth,
        xlabel={$q_1(t)$},
        ylabel={$q_2(t)$},
        ylabel style={yshift=-5pt},
        xmin=-1.5, xmax=1.0,
        ymin=-0.5, ymax=1,
        legend pos=north east,
        legend style={
        draw=none,           
        legend cell align=left}, 
        grid=both,
        major grid style={line width=0.4pt, opacity=0.4},
        minor grid style={line width=0.1pt, opacity=0.2},
        set layers
    ]
    \addplot+[very thick,
        color=black!80,
        only marks,
        mark size=2.5pt,
        mark=+
    ] table[
        header=true,
        x index=0,
        y index=1,
        col sep=comma
    ]{Data/motiv_ds.csv};
    \addlegendentry{Exact}
    \addplot+[very thick,
        color=blue!80,
        opacity=0.7,
        mark=none,
        mark size=2.5pt,
        skip coords between index={7}{100}, on layer=axis foreground
    ] table[
        header=true,
        x index=0,
        y index=1,
        col sep=comma
    ]{Data/motiv_pod.csv};
    \addlegendentry{POD}
    \addplot+[very thick,
        color=red!80,
        opacity=0.7,
        mark=none,
        mark size=2.5pt,
        skip coords between index={7}{100}, on layer=axis foreground
    ] table[
        header=true,
        x index=0,
        y index=1,
        col sep=comma
    ]{Data/motiv_ds.csv};
    \addlegendentry{Dynamic}[[ ]
    \addplot[-latex, red!80, ultra thick, quiver={u=0.94188722/1.5, v=-0.33592926/1.5}, on layer=axis foreground] coordinates {(0,0)};
    \node[draw, rectangle, fill=white, draw=none, on layer=axis foreground] (box2) at (0.5,-0.35){\color{red!80} \small $\mathbf{V}(\tau)$};
    \node[draw, rectangle, fill=white] (box2) at (-1.1,-0.35){\small $\tau=0.25$};
    \end{axis}
    \end{tikzpicture}
    \end{subfigure}
    \begin{subfigure}[b]{0.4\linewidth}
    \centering
    \begin{tikzpicture}
    \begin{axis}[
        width=\linewidth,
        xlabel={$q_1(t)$},
        ylabel={$q_2(t)$},
        ylabel style={yshift=-5pt},
        xmin=-1.5, xmax=1.0,
        ymin=-0.5, ymax=1,
        legend pos=north east,
        legend style={
        draw=none,           
        legend cell align=left}, 
        grid=both,
        major grid style={line width=0.4pt, opacity=0.4},
        minor grid style={line width=0.1pt, opacity=0.2},
        set layers
    ]
    \addplot+[very thick, color=black!80,
        only marks,
        mark size=2.5pt,
        mark=+
    ] table[
        header=true,
        x index=0,
        y index=1,
        col sep=comma
    ]{Data/motiv_ds.csv};
    \addlegendentry{Exact}
    \addplot+[very thick,
        opacity=0.7,
        color=blue!80,
        mark=none,
        mark size=2.5pt,
        skip coords between index={13}{100}, on layer=axis foreground
    ] table[
        header=true,
        x index=0,
        y index=1,
        col sep=comma
    ]{Data/motiv_pod.csv};
    \addlegendentry{POD}
    \addplot+[very thick,
        color=red!80,
        opacity=0.7,
        mark=none,
        mark size=2.5pt,
        skip coords between index={13}{100}, on layer=axis foreground
    ] table[
        header=true,
        x index=0,
        y index=1,
        col sep=comma
    ]{Data/motiv_ds.csv};
    \addlegendentry{Dynamic}
    \node[draw, rectangle, fill=white, draw=none] (box2) at (0.55,-0.25){\color{red!80} \small $\mathbf{V}(\tau)$};
    \node[draw, rectangle, fill=white] (box2) at (-1.1,-0.35) {\small $\tau=0.50$};
    \addplot[-latex, red!80, ultra thick, quiver={u=0.99382159/1.5, v=-0.11098942/1.5}, on layer=axis foreground] coordinates {(0,0)};
    \end{axis}
    \end{tikzpicture}
    \end{subfigure} \\[0.25em]
    \begin{subfigure}[b]{0.4\linewidth}
    \centering
    \begin{tikzpicture}
    \begin{axis}[
        width=\linewidth,
        xlabel={$q_1(t)$},
        ylabel={$q_2(t)$},
        ylabel style={yshift=-5pt},
        xmin=-1.5, xmax=1.0,
        ymin=-0.5, ymax=1,
        legend pos=north east,
        legend style={
        draw=none,           
        legend cell align=left}, 
        grid=both,
        major grid style={line width=0.4pt, opacity=0.4},
        minor grid style={line width=0.1pt, opacity=0.2},
        set layers
    ]
    \addplot[very thick, color=black!80,
        mark size=2.5pt,
        only marks,
        mark=+,  on layer=axis background
    ] table[
        header=true,
        x index=0,
        y index=1,
        col sep=comma
    ]{Data/motiv_ds.csv};
    \addlegendentry{Exact}
    \addplot+[very thick,
        opacity=0.7,
        mark size=2.5pt,
        color=blue!80,
        mark=none,
        skip coords between index={19}{100}, on layer=axis foreground
    ] table[
        header=true,
        x index=0,
        y index=1,
        col sep=comma
    ]{Data/motiv_pod.csv};
    \addlegendentry{POD}
    \addplot+[very thick,
        color=red!80,
        opacity=0.7,
        mark size=2.5pt,
        mark=none,
        skip coords between index={19}{100}, on layer=axis foreground
    ] table[
        header=true,
        x index=0,
        y index=1,
        col sep=comma
    ]{Data/motiv_ds.csv};
    \addlegendentry{Dynamic}
    \addplot[-latex, red!80, ultra thick, quiver={u=0.99279013/1.5, v=0.1198656/1.5}, on layer=axis foreground] coordinates {(0,0)};
    \node[draw, rectangle, fill=white, draw=none] (box) at (0.6,-0.11){\color{red!80} \small $\mathbf{V}(\tau)$};
    \node[draw, rectangle, fill=white] (box2) at (-1.1,-0.35) {\small $\tau=0.75$};
    \end{axis}
    \end{tikzpicture}
    \end{subfigure}
    \begin{subfigure}[b]{0.4\linewidth}
    \centering
    \begin{tikzpicture}
    \begin{axis}[
        width=\linewidth,
        xlabel={$q_1(t)$},
        ylabel={$q_2(t)$},
        ylabel style={yshift=-5pt},
        xmin=-1.5, xmax=1.0,
        ymin=-0.5, ymax=1,
        legend pos=north east,
        legend style={
        draw=none,           
        legend cell align=left}, 
        grid=both,
        major grid style={line width=0.4pt, opacity=0.4},
        minor grid style={line width=0.1pt, opacity=0.2},
        set layers
    ]
    \addplot+[very thick, color=black!80,
        mark size=2.5pt,
        only marks,
        mark=+] table[
        header=true,
        x index=0,
        y index=1,
        col sep=comma
    ]{Data/motiv_ds.csv};
    \addlegendentry{Exact}
    \addplot+[very thick,
        mark size=2.5pt,
        color=blue!80,
        opacity=0.7,
        mark=none, on layer=axis foreground
    ] table[
        header=true,
        x index=0,
        y index=1,
        col sep=comma
    ]{Data/motiv_pod.csv};
    \addlegendentry{POD}
    \addplot+[very thick,
        mark size=2.5pt,
        color=red!80,
        opacity=0.7,
        mark=none, on layer=axis foreground
    ] table[
        header=true,
        x index=0,
        y index=1,
        col sep=comma
    ]{Data/motiv_ds.csv};
    \addlegendentry{Dynamic}
    \addplot[-latex, red!80, ultra thick, quiver={u=0.93884781/1.5, v=0.34433238/1.5}, on layer=axis foreground] coordinates {(0,0)};
    \node[draw, rectangle, fill=white, draw=none] (box2) at (0.29,0.29){\color{red!80} \small $\mathbf{V}(\tau)$};
    \node[draw, rectangle, fill=white] (box2) at (-1.1,-0.35) {\small $\tau=1.00$};
    \end{axis}
    \end{tikzpicture}
    \end{subfigure}
    \caption{Comparing a time-continuous subspace approximation of the form~\cref{eq:model_form} (red) and a static POD approximation of the form~\cref{eq:pod} (blue) for a 2D trajectory at different time steps, illustrating the framework's ability to approximate solution-manifold curvature. The direction of $\mathbf{V}(\tau)$ is represented by the red arrow.}
\label{fig:motivational}
\end{figure}

We now demonstrate the merit of time-continuous subspace approximation using a two-dimensional example. The system dimension $N = 2$ is already very small, but a reduction to $r = 1$ dimensions provides valuable geometric intuition. Consider the trajectory
\begin{equation}
    \mathbf{q}(t) = \left[\begin{array}{c}
        q_1(t) \\[2pt] q_2(t)
    \end{array}\right]
    = \left[\begin{array}{c}
        t \\[2pt] \dfrac{t^2}{2}
    \end{array}\right],
    \qquad
    t \in [-1.25,0.75].
\end{equation}
The reference state is chosen to be the origin, $\bar{\mathbf{q}} = [~0~~0~]^\top$, which lies on the trajectory $\mathbf{q}(t)$ at $t = 0$.
From a dataset of $K = 25$ snapshots uniformly spaced in time, POD finds the first principal component of the data (blue curve in \Cref{fig:motivational}) to be the vector $[~-0.937~~0.350~]^\top$.
The resulting POD state approximation is therefore
\begin{equation}
    \mathbf{q}(t) \approx 
    \underbrace{\left[\begin{array}{r}
        -0.937 \\
        0.350
    \end{array}\right]}_{\textstyle \mathbf{V}}
    \hat{q}(t),
\end{equation}
where $\hat{q}(t)$ is the \emph{scalar} temporal coefficient at time $t$.

In contrast, a dynamic approximation \cref{eq:dynamic} changes the basis with time in a way that tracks the  solution manifold. This basis rotation follows the shortest possible path between subspaces, ensuring the representation evolves along the manifold's natural curvature. By allowing the basis to rotate and align itself with the instantaneous direction of the system state, the model effectively decouples the geometric rotation of the coordinate system from the actual coordinate scaling. This behavior is visualized in \Cref{fig:motivational}.

Setting the number of features to $n_f=2$ allows the dynamic approximation to span $\mathbb{R}^2$. The Riemannian optimization problem~\cref{eq:riemannian_optimization} converges with a tolerance of $10^{-4}$ in $16$ optimization iterations, learning an approximation of the form
\begin{equation}
    \mathbf{q}(t) \approx \underbrace{\underbrace{\left[\begin{array}{rr}
        -0.974 & 0.226 \\
         0.226 & 0.974
    \end{array}\right]}_{\textstyle \boldsymbol{\Omega}}
    \ \underbrace{\left[\begin{array}{c}
    \cos(2.562\tau) \\
    \sin(2.562\tau)
    \end{array}\right]}_{\textstyle \mathbf{W}(\tau,\boldsymbol{\Theta})}}_{\textstyle \mathbf{V}(\tau)} \hat{q}(t).
\end{equation}
The coefficient $\hat{q}(t)$ at any given time $t$ may be computed by projecting the snapshot data onto the instantaneous basis through $\hat{q}(t) = \mathbf{V}(\tau)^\top \mathbf{q}(t)$. The relative error for the POD approximation was found to equal $27.8\%$, whereas the error reduces to $1.1\%$ for the dynamic subspace approximation. 

Clearly, the dynamic subspace approach effectively learns and represents curved solution manifolds by constructing time-continuous approximations. While this particular model problem can certainly be solved with established nonlinear approximation methods, such as quadratic manifolds~\cite{barnett2022quadratic, geelen2023operator, geelen2024learning} or autoencoders~\cite{lee2020model}, the authors contend that the proposed framework offers a unique and computationally efficient alternative that directly exploits the underlying dynamical structure.

\subsection{Overview}

\begin{algorithm}[!tbp]
\caption{Scalable dynamic subspace estimation.}
\label{alg:dynamic_subspace}
\begin{algorithmic}[1]
    \Require Snapshot matrix $\mathbf{Q} \in \mathbb{R}^{N \times K}$; reference state $\bar{\mathbf{q}}\in\mathbb{R}^{N}$; target rank $r$; number of features $n_f$

    \LineComment{Step 1: Feature extraction}
    \State $\mathbf{Q}' \gets \mathbf{Q} - \bar{\mathbf{q}}\mathbf{1}^\top$
        \Comment{Center snapshots by the reference state.}
    \State $\mathbf{X}\boldsymbol{\Xi}\mathbf{Z}^\top \gets \operatorname{SVD}(\mathbf{Q}')$
        \Comment{Compute thin SVD of centered snapshot data.}
    \State $\tilde{\mathbf{V}} \gets \mathbf{X}_{:,:n_f}\in\operatorname{St}(N, n_f)$ 
        \Comment{Extract leading $n_f$ left singular vectors.}
    \State $\mathbf{Y} \gets \tilde{\mathbf{V}}^\top \mathbf{Q}' \in \mathbb{R}^{n_f \times K}$
        \Comment{Project snapshots into feature space.}

    \LineComment{Step 2: Riemannian optimization}
    \State $\mathbf{T} \gets \text{random sample from}~\operatorname{St}(n_f, 2r)$
        \Comment{Initial guess for orthonormal matrices.}
    \State $\boldsymbol{\Theta} \gets \text{random sample from}~\mathbb{R}^r$
        \Comment{Initial guess for principal angles.}
    \State $\mathcal{M} \gets \operatorname{St}(n_f, 2r) \times \mathbb{R}^r$
        \Comment{Optimization search space (product manifold).}
    \State Minimize the normalized reconstruction error in feature space via \texttt{Pymanopt}~\cite{townsend2016pymanopt} with the \texttt{JAX}~\cite{jax2018github} backend:
    \begin{equation*}
        \mathbf{T}_\text{opt}, \boldsymbol{\Theta}_\text{opt}
        \gets \argmin_{( \mathbf{T}, \boldsymbol{\Theta}) \in \mathcal{M}} 
        \frac{\|\mathbf{Y} - \mathbf{T}\mathbf{R}(\mathbf{T}, \boldsymbol{\Theta}, \mathbf{Y})\|_F}{\|\mathbf{Y}\|_F},
    \end{equation*}
    \item[] where $\mathbf{R}(\mathbf{T},\boldsymbol{\Theta}, \mathbf{Y}) \in \mathbb{R}^{2r \times K}$, defined in \cref{eq:R}, is updated at every iteration of the optimizer.

    \LineComment{Step 3: Assemble the dynamic basis}
    \State $\mathbf{V}(\tau) \gets \tilde{\mathbf{V}} \mathbf{T}_\text{opt} [~\cos(\tau\boldsymbol{\Theta}_\text{opt})~~\sin(\tau\boldsymbol{\Theta}_\text{opt})~]\in\operatorname{St}(N, r)$
    \vspace{.5em}
    \Ensure Dynamic basis matrix $\mathbf{V}(\tau)$
\end{algorithmic}
\end{algorithm}

The proposed dynamic subspace approach is implemented through a scalable Riemannian optimization framework that decouples the computational cost from the high-dimensional state space $N$. To ensure numerical efficiency, the optimization is performed within a reduced feature space spanned by $n_f$ dominant modes extracted via an singular value decomposition of the snapshot data. The algorithm treats the search space as a product manifold $\mathcal{M} = \operatorname{St}(n_f, 2r) \times \mathbb{R}^r$, consisting of the Stiefel manifold of orthonormal matrices and the Euclidean space of principal angles. Specifically, the algorithm accepts the high-dimensional snapshot data $\mathbf{Q}$, a reference state $\bar{\mathbf{q}}$, the desired reduced rank $r$, and the number of features $n_f$ as inputs; it then identifies a set of optimal dynamic basis vectors that vary continuously with time according to the normalized time $\tau$. Using the fully vectorized formulation introduced in \Cref{subsec:vectorization}, the framework computes the projections for all snapshots in a single pass, effectively eliminating the overhead associated with temporal looping. The training phase leverages a Riemannian conjugate gradient solver to minimize the normalized Frobenius norm of the reconstruction error, employing retraction mappings and vector transport to ensure the evolving basis remains strictly orthonormal while respecting the natural curvature of the Grassmannian geodesic. The specific computational steps of this procedure are summarized in \Cref{alg:dynamic_subspace}.

\section{Results}
\label{sec:results}

To evaluate the performance and scalability of the proposed dynamic subspace framework, we conduct numerical experiments on two distinct problems that characterize the challenges of the Kolmogorov barrier. The first experiment focuses on a 1D transport-dominated problem, which serves as a benchmark for assessing the ability of the geodesic-based approach to track traveling waves without the rank inflation required by static POD. The second experiment demonstrates the method's scalability on a high-fidelity turbulent airfoil wake dataset. This large-scale case involves complex, non-stationary physics, highlighting the efficiency of the feature-space optimization and the vectorized implementation. In both cases, we compare the dynamic basis against traditional static POD to quantify gains in both reconstruction accuracy and coordinate sparsity. 

\subsection{One-dimensional Transport Equation}
\label{subsec:one_d_transport}

The one-dimensional linear transport equation is a fundamental PDE that describes how a scalar quantity is carried along by a fluid moving at a constant velocity.
Mathematically expressed as 
\begin{equation}
    \dfrac{\partial s}{\partial t} + c \dfrac{\partial s}{\partial x} = 0,
\end{equation}
it represents a conservative system where an initial profile remains unchanged in shape as it translates across space. Setting the initial condition to a Gaussian pulse centered at $\mu\in\mathbb{R}$ with width parameter $\sigma$,
\begin{equation}
    s(x,t_0)
    = s_0(x)
    = \dfrac{1}{\sigma\sqrt{2\pi}} \exp\left( - \dfrac{(x - \mu)^2}{2\sigma^2} \right),
    \quad x \in [0, 1],
\end{equation}
the analytic solution is given by $s(x, t) = s_0(x - ct)$.
As time progresses, the pulse advects positively (away from $x = \mu$) with wave speed $c$.

\begin{figure}[!t]
    \centering \scriptsize
    \begin{subfigure}{\linewidth}
    \begin{subfigure}[b]{0.49\linewidth}
    \hspace{-1em}
  \begin{tikzpicture}
    \begin{axis}[
        width=\linewidth,
        height=.5\linewidth,
        xlabel={$x$-coordinate},
        ylabel={Mode shape},
        ylabel style={yshift=-5pt},
        xlabel style={xshift=-5pt},
        xmin=0, xmax=1.0,
        ymin=-0.12, ymax=0.12,
        ytick={-0.15,-0.1,...,0.15},
        yticklabel style={
        /pgf/number format/fixed,
        /pgf/number format/precision=2,
        /pgf/number format/fixed zerofill 
        },
        scaled y ticks=false,
        legend style={
        at={(0.5, 1.01)},      
        anchor=south,          
        legend columns=-1,     
        draw=none,             
        /tikz/every even column/.append style={column sep=1pt}, 
        draw=none,           
        legend cell align=left},
        grid=both,
        major grid style={line width=0.4pt, opacity=0.4},
        minor grid style={line width=0.1pt, opacity=0.2},
        set layers,
        colormap={example}{
        samples of colormap={
            5 of viridis,
        },
        },
        legend image code/.code={
        \draw[mark repeat=2, mark phase=2]
        plot coordinates {
            (0cm,0cm)
            (0.3cm,0cm) 
        };
    }
    ]

\addplot[thick, mark=none, index of colormap=0] table[header=true, x index=0, y index=1, col sep=comma]{Data/dynamic_vectors_t0.csv};
\addlegendentry{$t=0$}

\addplot[thick, mark=none, index of colormap=1] table[header=true, x index=0, y index=1, col sep=comma]{Data/dynamic_vectors_t1.csv};
\addlegendentry{$t=0.02$}

\addplot[thick, mark=none, index of colormap=2] table[header=true, x index=0, y index=1, col sep=comma]{Data/dynamic_vectors_t2.csv};
\addlegendentry{$t=0.04$}

\addplot[thick, mark=none, index of colormap=4] table[header=true, x index=0, y index=1, col sep=comma]{Data/dynamic_vectors_t3.csv};
\addlegendentry{$t=0.06$}

\end{axis}
\end{tikzpicture}
    \end{subfigure}
    \begin{subfigure}[b]{0.49\linewidth}
    \hspace{-1em}
    \begin{tikzpicture}
    \begin{axis}[
        width=\linewidth,
        height=.5\linewidth,
        xlabel={$x$-coordinate},
        ylabel={Mode shape},
        ylabel style={yshift=-5pt},
        xlabel style={xshift=-5pt},
        xmin=0, xmax=1.0,
        ymin=-0.12, ymax=0.12,
        ymin=-0.12, ymax=0.12,
        ytick={-0.15,-0.1,...,0.15},
        yticklabel style={
        /pgf/number format/fixed,
        /pgf/number format/precision=2,
        /pgf/number format/fixed zerofill 
        },
        scaled y ticks=false,
        grid=both,
        major grid style={line width=0.4pt, opacity=0.4},
        minor grid style={line width=0.1pt, opacity=0.2},
        set layers
    ]
    \addplot+[thick,
        color=black!80,
        mark=none
    ] table[
        header=true,
        x index=0,
        y index=1,
        col sep=comma
    ]{Data/transport_pod.csv};
    \end{axis}
    \end{tikzpicture}
    \end{subfigure}
    \vspace{-0.5em}
    \caption{Dynamic mode \#1 (left) and the dominant POD mode (right).}
    \end{subfigure}\\[0.5em]
    \begin{subfigure}{\linewidth}
    \begin{subfigure}[b]{0.49\linewidth}
    \hspace{-1em}
    \begin{tikzpicture}
    \begin{axis}[
        width=\linewidth,
        height=.5\linewidth,
        xlabel={$x$-coordinate},
        ylabel={Mode shape},
        ylabel style={yshift=-5pt},
        xlabel style={xshift=-5pt},
        xmin=0, xmax=1.0,
        ymin=-0.12, ymax=0.12,
        ytick={-0.15,-0.1,...,0.15},
        yticklabel style={
        /pgf/number format/fixed,
        /pgf/number format/precision=2,
        /pgf/number format/fixed zerofill 
        },
        scaled y ticks=false,
        legend style={
        at={(0.5, 1.01)},      
        anchor=south,          
        legend columns=-1,     
        draw=none,             
        /tikz/every even column/.append style={column sep=1pt}, 
        draw=none,           
        legend cell align=left},
        grid=both,
        major grid style={line width=0.4pt, opacity=0.4},
        minor grid style={line width=0.1pt, opacity=0.2},
        set layers,
        colormap={example}{
        samples of colormap={
            5 of viridis,
        },
        },
        legend image code/.code={
        \draw[mark repeat=2, mark phase=2]
        plot coordinates {
            (0cm,0cm)
            (0.3cm,0cm) 
        };
    }
    ]
\addplot[thick, mark=none, index of colormap=0] table[header=true, x index=0, y index=10, col sep=comma]{Data/dynamic_vectors_t0.csv};
\addlegendentry{$t=0$}

\addplot[thick, mark=none, index of colormap=1] table[header=true, x index=0, y index=10, col sep=comma]{Data/dynamic_vectors_t1.csv};
\addlegendentry{$t=0.02$}

\addplot[thick, mark=none, index of colormap=2] table[header=true, x index=0, y index=10, col sep=comma]{Data/dynamic_vectors_t2.csv};
\addlegendentry{$t=0.04$}

\addplot[thick, mark=none, index of colormap=4] table[header=true, x index=0, y index=10, col sep=comma]{Data/dynamic_vectors_t3.csv};
\addlegendentry{$t=0.06$}

\end{axis}
    \end{tikzpicture}
    \end{subfigure}
    \begin{subfigure}[b]{0.49\linewidth}
    \hspace{-1em}
    \begin{tikzpicture}
    \begin{axis}[
        width=\linewidth,
        height=.5\linewidth,
        xlabel={$x$-coordinate},
        ylabel={Mode shape},
        ylabel style={yshift=-5pt},
        xlabel style={xshift=-5pt},
        xmin=0, xmax=1.0,
        ymin=-0.12, ymax=0.12,
        ytick={-0.15,-0.1,...,0.15},
        yticklabel style={
        /pgf/number format/fixed,
        /pgf/number format/precision=2,
        /pgf/number format/fixed zerofill 
        },
        scaled y ticks=false,
        legend pos=north east,
        grid=both,
        major grid style={line width=0.4pt, opacity=0.4},
        minor grid style={line width=0.1pt, opacity=0.2},
        set layers,
    ]
    \addplot+[thick,
        color=black!80,
        mark=none
    ] table[
        header=true,
        x index=0,
        y index=10,
        col sep=comma
    ]{Data/transport_pod.csv};
    \end{axis}
    \end{tikzpicture}
    \end{subfigure}
    \vspace{-0.5em}
    \caption{Dynamic mode \#10 (left) and and POD mode \#10 (right).}
    \end{subfigure}\\[0.5em]
    \begin{subfigure}{\linewidth}
    \begin{subfigure}[b]{0.49\linewidth}
    \hspace{-1em}
    \begin{tikzpicture}
    \begin{axis}[
        width=\linewidth,
        height=.5\linewidth,
        xlabel={$x$-coordinate},
        ylabel={Mode shape},
        ylabel style={yshift=-5pt},
        xlabel style={xshift=-5pt},
        xmin=0, xmax=1.0,
        ymin=-0.12, ymax=0.12,
        ytick={-0.15,-0.1,...,0.15},
        yticklabel style={
        /pgf/number format/fixed,
        /pgf/number format/precision=2,
        /pgf/number format/fixed zerofill 
        },
        scaled y ticks=false,
        legend style={
        at={(0.5, 1.01)},      
        anchor=south,          
        legend columns=-1,     
        draw=none,             
        /tikz/every even column/.append style={column sep=1pt}, 
        draw=none,           
        legend cell align=left},
        grid=both,
        major grid style={line width=0.4pt, opacity=0.4},
        minor grid style={line width=0.1pt, opacity=0.2},
        set layers,
        colormap={example}{
        samples of colormap={
            5 of viridis,
        },
        },
        legend image code/.code={
        \draw[mark repeat=2, mark phase=2]
        plot coordinates {
            (0cm,0cm)
            (0.3cm,0cm) 
        };
    }
    ]
\addplot[thick, mark=none, index of colormap=0] table[header=true, x index=0, y index=20, col sep=comma]{Data/dynamic_vectors_t0.csv};
\addlegendentry{$t=0$}

\addplot[thick, mark=none, index of colormap=1] table[header=true, x index=0, y index=20, col sep=comma]{Data/dynamic_vectors_t1.csv};
\addlegendentry{$t=0.02$}

\addplot[thick, mark=none, index of colormap=2] table[header=true, x index=0, y index=20, col sep=comma]{Data/dynamic_vectors_t2.csv};
\addlegendentry{$t=0.04$}

\addplot[thick, mark=none, index of colormap=4] table[header=true, x index=0, y index=20, col sep=comma]{Data/dynamic_vectors_t3.csv};
\addlegendentry{$t=0.06$}

    \end{axis}
    \end{tikzpicture}
    \end{subfigure}
    \begin{subfigure}[b]{0.49\linewidth}
    \hspace{-1em}
    \begin{tikzpicture}
    \begin{axis}[
        width=\linewidth,
        height=.5\linewidth,
        xlabel={$x$-coordinate},
        ylabel={Mode shape},
        ylabel style={yshift=-5pt},
        xlabel style={xshift=-5pt},
        xmin=0, xmax=1.0,
        ymin=-0.12, ymax=0.12,
        ytick={-0.15,-0.1,...,0.15},
        yticklabel style={
        /pgf/number format/fixed,
        /pgf/number format/precision=2,
        /pgf/number format/fixed zerofill 
        },
        scaled y ticks=false,
        grid=both,
        major grid style={line width=0.4pt, opacity=0.4},
        minor grid style={line width=0.1pt, opacity=0.2},
        set layers
    ]
    \addplot+[thick,
        color=black!80,
        mark=none
    ] table[
        header=true,
        x index=0,
        y index=20,
        col sep=comma
    ]{Data/transport_pod.csv};
    \end{axis}
    \end{tikzpicture}
    \end{subfigure}
    \vspace{-0.5em}
    \caption{Dynamic mode \#20 (left) and POD mode \#20 (right).}
    \end{subfigure}
    \caption{Temporal evolution of the selected dynamic basis at a dimensionality of $r=20$. The static POD modes are shown for reference.}
    \label{fig:dynamic_modes_1d}
\end{figure}

\begin{figure}[!t]
    \centering \small
    \begin{tikzpicture}
    \begin{axis}[
        width=.65\linewidth,
        height=.5*.65\linewidth,
        xlabel={$x$-coordinate},
        ylabel={Mode shape},
        xmin=0, xmax=1.0,
        ymin=-5, ymax=30,
        scaled y ticks=false,
        legend style={
        at={(0.5, 1.01)},      
        anchor=south,          
        legend columns=-1,     
        draw=none,             
        /tikz/every even column/.append style={column sep=1pt}, 
        draw=none,           
        legend cell align=left},
        grid=both,
        major grid style={line width=0.4pt, opacity=0.4},
        minor grid style={line width=0.1pt, opacity=0.2},
        set layers,
        colormap={example}{
        samples of colormap={
            4 of viridis,
        },
        }
        ]
    \addplot[very thick, mark=+, only marks, index of colormap=0, mark repeat=6] table[header=true, x index=0, y index=1, col sep=comma]{Data/transport_exact_state.csv};
    \addlegendentry{Exact}
    \addplot[very thick, mark=+, only marks, index of colormap=1, mark repeat=6, forget plot] table[header=true, x index=0, y index=2, col sep=comma]{Data/transport_exact_state.csv};
    \addplot[very thick, mark=+, only marks, index of colormap=2, mark repeat=6, forget plot] table[header=true, x index=0, y index=3, col sep=comma]{Data/transport_exact_state.csv};
    \addplot[very thick, mark=+, only marks, index of colormap=3, mark repeat=6, forget plot] table[header=true, x index=0, y index=4, col sep=comma]{Data/transport_exact_state.csv};
    \addplot[very thick, mark=none, index of colormap=0, densely dashed] table[header=true, x    
    index=0, y index=1, col sep=comma]{Data/transport_pod_state.csv};
    \addlegendentry{POD}
    \addplot[very thick, mark=none, index of colormap=1, densely dashed, forget plot] table[header=true, x index=0, y index=2, col sep=comma]{Data/transport_pod_state.csv};
    \addplot[very thick, mark=none, index of colormap=2, densely dashed, forget plot] table[header=true, x index=0, y index=3, col sep=comma]{Data/transport_pod_state.csv};
    \addplot[very thick, mark=none, index of colormap=3, densely dashed, forget plot] table[header=true, x index=0, y index=4, col sep=comma]{Data/transport_pod_state.csv};
    \addplot[very thick, mark=none, index of colormap=0] table[header=true, x    
    index=0, y index=1, col sep=comma]{Data/transport_ds_state.csv};
    \addlegendentry{Dynamic}
    \addplot[very thick, mark=none, index of colormap=1, forget plot] table[header=true, x index=0, y index=2, col sep=comma]{Data/transport_ds_state.csv};
    \addplot[very thick, mark=none, index of colormap=2, forget plot] table[header=true, x index=0, y index=3, col sep=comma]{Data/transport_ds_state.csv};
    \addplot[very thick, mark=none, index of colormap=3, forget plot] table[header=true, x index=0, y index=4, col sep=comma]{Data/transport_ds_state.csv};
    \node[fill=white] (box1) at (0.17,20){\footnotesize $t=0$};
    \node[fill=white] (box1) at (0.435,20 ){\footnotesize $t=0.02$};
    \node[fill=white] (box1) at (0.685,20){\footnotesize $t=0.04$};
    \node[fill=white] (box1) at (0.93,20) {\footnotesize $t=0.06$};
    \end{axis}
    \end{tikzpicture}
    \vspace{-0.5em}
    \caption{Qualitative comparison between a dynamic subspace and a (static) POD subspace for reconstructing time-series data at $r=20$ degrees of freedom.}
    \label{fig:comparison_1d}
\end{figure}
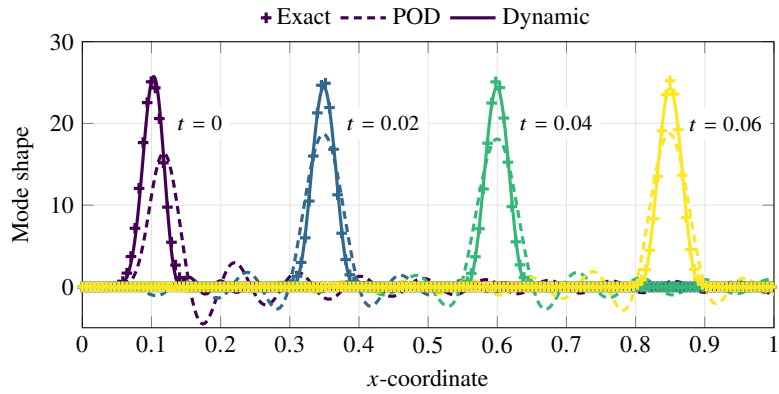

For a numerical experiment, we consider the time domain $t \in [0, 0.1]$, the spatial domain $x\in[0,1]$, and set the wave speed to $c = 10.0$, the center of the initial condition to $\mu = 0.1$, and the initial condition width to $\sigma = 0.01$.
We then evaluate the function $s(x,t)$ at temporal instances $t_j = 0.2(j - 1)/500$ for $j = 1, \dots, 500$ over $N=1{,}024$ equidistant spatial points $x \in [0, 1]$ and organize the results into training and testing data matrices,
\begin{equation}
    \mathbf{Q}_{\text{train}} = [~\mathbf{q}(t_{1})~~\mathbf{q}(t_{3})~~\cdots~~\mathbf{q}(t_{499})~] \in \mathbb{R}^{1024 \times 250},
    \qquad
    \mathbf{Q}_{\text{test}} = [~\mathbf{q}(t_{2})~~\mathbf{q}(t_{4})~~\cdots~~\mathbf{q}(t_{500})~] \in \mathbb{R}^{1024 \times 250}.    
\end{equation}
The Riemannian optimization employs a stopping criterion based on a minimum gradient norm of $10^{-4}$, with the initial guess determined by a random set of principal angles. To ensure numerical robustness and consistency, the approximation is learned five times using different initializations, and the average error across these runs is reported. The reference state for this dataset is set to the zero vector, $\bar{\mathbf{q}} = \mathbf{0}$, and \Cref{alg:dynamic_subspace} is applied using $\mathbf{Q} = \mathbf{Q}_\text{train}$, $n_f = 200$, and for several values of $r$.

\Cref{fig:dynamic_modes_1d} illustrates the temporal evolution of the dynamic basis vectors for a reduced dimension of $r=20$ alongside the static POD modes for comparative analysis. While POD modes follow a distinct hierarchy in which higher-order modes progressively capture finer features, this rigid ordering is absent in the proposed geodesic approach. Instead, fine-scale structures are distributed across the dynamic basis vectors as they adapt to the underlying geometry of the manifold. This flexibility allows the representation to maintain high fidelity without adhering to the energy-based ranking characteristic of POD. As shown in \Cref{fig:comparison_1d}, the proposed approach offers a clear qualitative advantage at $r=20$ degrees of freedom. It effectively captures the amplitude of the moving peak while simultaneously suppressing oscillations in its wake, a performance gain that remains consistent across the entire observation window.

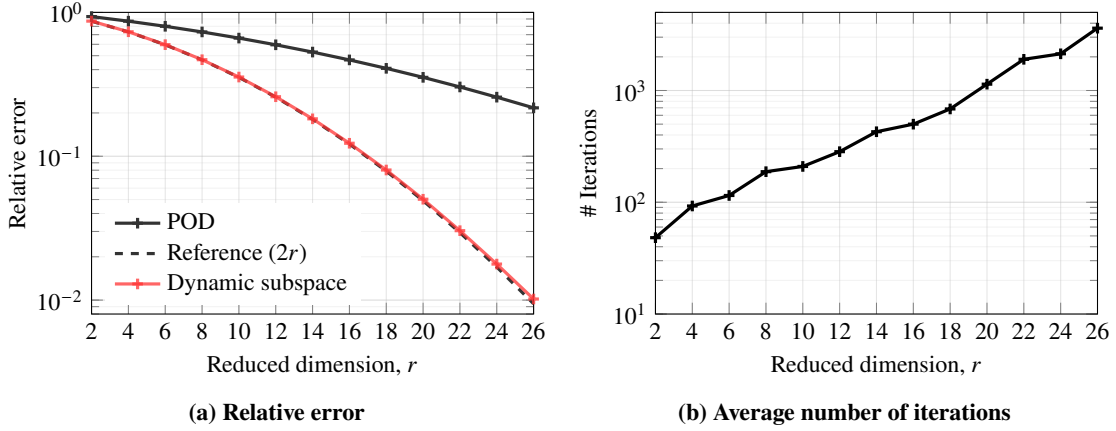
\begin{figure}[!tbp]
\centering \small
\begin{subfigure}[t]{0.45\textwidth}
    \centering
    \begin{tikzpicture}
    \begin{semilogyaxis}[
        width=\linewidth,
        height=.75\linewidth,
        xlabel={Reduced dimension, $r$},
        ylabel={Relative error},
        ylabel style={yshift=-2pt},
        xmin=2, xmax=26,
        ymin=0.008, ymax=1,
        xtick={2,4,...,26},
        legend pos=south west,
        legend style={
        draw=none,           
        legend cell align=left}, 
        grid=both,
        major grid style={line width=0.4pt, opacity=0.4},
        minor grid style={line width=0.1pt, opacity=0.2},
    ]
    \addplot[very thick,
        color=black!80,
        mark=+,
    ] table[
        x index=0,
        y index=1,
        col sep=comma
    ]{Data/1d_results.csv};
    \addlegendentry{POD}
    \addplot+[dashed, very thick,
        color=black!80,
        mark=none,
    ] table[
        x index=0,
        y index=2,
        col sep=comma
    ]{Data/1d_results.csv};
    \addlegendentry{Reference ($2r$)}
    \addplot+[very thick,
        color=red!80,
        draw opacity=0.75,
        mark=+,
    ] table[
        header=true,
        x index=0,
        y index=8,
        col sep=comma
    ]{Data/1d_results.csv};
    \addlegendentry{Dynamic subspace}
    \end{semilogyaxis}
    \end{tikzpicture}
    \caption{Relative error}
    \label{fig:one_dimensional_error}
\end{subfigure}
\begin{subfigure}[t]{0.45\textwidth}
    \centering
    \begin{tikzpicture}
    \begin{semilogyaxis}[
        width=\linewidth,
        height=.75\linewidth,
        xlabel={Reduced dimension, $r$},
        ylabel={\# Iterations},
        xmin=2, xmax=26,
        ymin=10, ymax=5e3,
        xtick={2,4,...,26},
        mark size=2pt,
        grid=both,
        major grid style={line width=0.4pt, opacity=0.4},
        minor grid style={line width=0.1pt, opacity=0.2},
    ]
    \addplot+[very thick,
        color=black,
        mark=+,
    ] table[
        x index=0,
        y index=14,
        col sep=comma
    ]{Data/1d_results.csv};
    \end{semilogyaxis}
    \end{tikzpicture}
    \caption{Average number of iterations}
    \label{fig:one_dimensional_iterations}
\label{nonlinear_dimension_comp}
\end{subfigure}
\caption{Relative error (left) and the average number of iterations required for the gradient-based optimization algorithm to converge to within a tolerance of $10^{-4}$ (right) for the one-dimensional transport problem.}
\label{relative_err}
\end{figure}

A quantitative comparison between the POD and the proposed approach is presented in \Cref{fig:one_dimensional_error}.
The representation error at each dimension $r$ in the proposed dynamic approach is computed as
\begin{equation}
     \frac{\|\mathbf{Y}_\text{test} - \mathbf{T}_\text{opt} \mathbf{R}(\mathbf{T}_\text{opt}, \boldsymbol{\Theta}_\text{opt}, \mathbf{Y}_\text{test}) \|_F}{\| \mathbf{Y}_\text{test} \|_F},
\end{equation}
in which $\mathbf{Y}_\text{test} = \tilde{\mathbf{V}}^\top\mathbf{Q}_\text{test}$ and $\tilde{\mathbf{V}}$ is computed in \Cref{alg:dynamic_subspace} from $\mathbf{Q}_\text{train}$.
For POD, the equivalent error metric is the relative projection error,
\begin{equation}
    \frac{\|\mathbf{Q}_\text{test} - \mathbf{ZZ}^\top\mathbf{Q}_\text{test}\|_F}{\|\mathbf{Q}_\text{test}\|_F},
\end{equation}
where $\mathbf{Z} \in \mathbb{R}^{N\times r}$ is the POD basis matrix computed from $\mathbf{Q}_\text{train}$.
For transport-dominated phenomena like the wave equation, it is well understood why the POD relative error converges slowly as dimensionality increases~\cite{peherstorfer2022breaking}. The proposed approach significantly mitigates this issue and demonstrates much faster convergence. Following the observation in~\cite{blocker2023dynamicsubspaceestimationgrassmannian} that a rank-$r$ geodesic spans a space with a dimension as large as $2r$, we include the subspace error for a $2r$-rank POD as a reference. The proposed method follows this idealized scenario very closely, validating the efficacy of the geodesic perspective. Finally, \Cref{fig:one_dimensional_iterations} reports the average number of iterations required for convergence. Although the optimization becomes more challenging as the number of variables grows, the computational complexity remains independent of the original state dimension $N$. This independence makes the approach fast and therefore highly suitable for nonlinear dimensionality reduction in large-scale engineering applications, as demonstrated in the subsequent experiment.

\subsection{Turbulent Airfoil Wake LES}

In this section, we examine a data set from the large eddy simulation database of~\cite{deepblue, towne2023database, yeh2019resolvent}, which provides a high‑fidelity three‑dimensional LES flow field.
This dataset consists of a turbulent separated flow at Mach 0.3 over a NACA 0012 airfoil at \ang{6} angle of attack. The Reynolds number of the flow is $23{,}000$, featuring a transitional boundary layer, separation over a recirculation bubble, and a turbulent wake.

\begin{figure}[!t]
    \centering \footnotesize
    \begin{subfigure}[t]{.58\linewidth}
        \includegraphics[width=\linewidth]{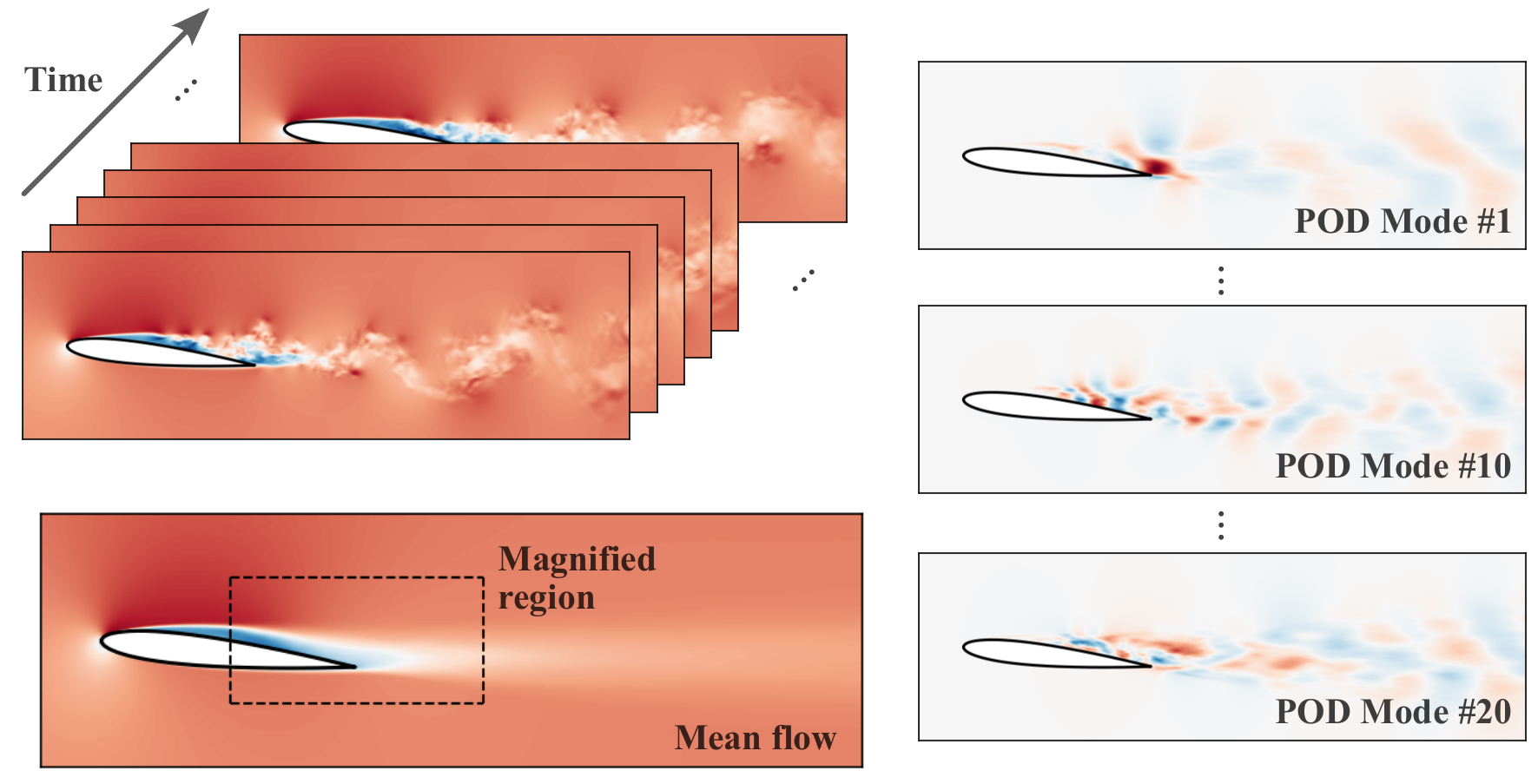}
        \caption{Modal decomposition of the fluid data}
    \end{subfigure}
    \hfill
    \begin{subfigure}[t]{.41\linewidth}
    \begin{tikzpicture}
    \begin{semilogyaxis}[
        width=\linewidth,
        height=.75\linewidth,
        xlabel={Index},
        ylabel={Normalized singular value},
        ylabel style={yshift=-2pt},
        xmin=0, xmax=800,
        ymin=0.01, ymax=1,
        mark repeat=10, 
        xtick={0,100,...,800},
        legend pos=south west,
        legend style={
        draw=none,           
        legend cell align=left}, 
        grid=both,
        major grid style={line width=0.4pt, opacity=0.4},
        minor grid style={line width=0.1pt, opacity=0.2},
    ]
    \addplot+[very thick,
        color=black!80,
        mark=+,
    ] table[
        header=true,
        x index=0,
        y index=2,
        col sep=comma
    ]{Data/les_sv.csv};
    \end{semilogyaxis}
    \end{tikzpicture}
    \caption{Singular values}
    \end{subfigure}
    \caption{Overview of the turbulent airfoil wake LES dataset from~\cite{deepblue}. Each snapshot provides all three velocity components extracted over the midspan slice. We plot the $x$-component of the velocity fields for some snapshots alongside the mean flow (left), a few POD modes of the training data (middle), and the singular values (right).}
    \label{fig:les_overview}
\end{figure}

\begin{figure}[!tbp]
\centering \small
\begin{subfigure}[t]{0.45\textwidth}
    \centering
    \begin{tikzpicture}
    \begin{semilogyaxis}[
        width=\linewidth,
        height=.75\linewidth,
        xlabel={Reduced dimension, $r$},
        ylabel={Relative error},
        ylabel style={yshift=-2pt},
        xmin=5, xmax=80,
        ymin=0.035, ymax=0.15,
        xtick={0,10,...,80},
        legend pos=south west,
        legend style={
        draw=none,           
        legend cell align=left}, 
        grid=both,
        major grid style={line width=0.4pt, opacity=0.4},
        minor grid style={line width=0.1pt, opacity=0.2},
    ]
    \addplot[very thick,
        color=black!80,
        mark=+,
    ] table[
        x index=0,
        y index=1,
        col sep=comma
    ]{Data/les_results.csv};
    \addlegendentry{POD}
    \addplot+[dashed, very thick,
        color=black!80,
        mark=none,
    ] table[
        x index=0,
        y index=2,
        col sep=comma
    ]{Data/les_results.csv};
    \addlegendentry{Reference ($2r$)}
    \addplot+[very thick,
        color=red!80,
        draw opacity=0.75,
        mark=+,
    ] table[
        header=true,
        x index=0,
        y index=8,
        col sep=comma
    ]{Data/les_results.csv};
    \addlegendentry{Dynamic subspace}
    \end{semilogyaxis}
    \end{tikzpicture}
    \caption{Relative error}
    \label{fig:les_error}
\end{subfigure}
\begin{subfigure}[t]{0.45\textwidth}
    \centering
    \begin{tikzpicture}
    \begin{semilogyaxis}[
        width=\linewidth,
        height=.75\linewidth,
        xlabel={Reduced dimension, $r$},
        ylabel={\# Iterations},
        xmin=5, xmax=80,
        ymin=10, ymax=1e3,
        xtick={0,10,...,80},
        mark size=2pt,
        grid=both,
        major grid style={line width=0.4pt, opacity=0.4},
        minor grid style={line width=0.1pt, opacity=0.2},
    ]
    \addplot+[very thick,
        color=black,
        mark=+,
    ] table[
        x index=0,
        y index=14,
        col sep=comma
    ]{Data/les_results.csv};
    \node[draw, rectangle, fill=white, draw=none] (box) at (65,350) {$\approx 25\ $sec.};
    \coordinate (target) at (78, 200);
    \draw[-latex, thick] (box) -- (target) node[midway, above] {};
    \node[draw, rectangle, fill=white, draw=none] (box2) at (22,35) {$\approx 0.24\ $sec.};
    \coordinate (target2) at (7, 24);
    \draw[-latex, thick] (box2) -- (target2) node[near end, right] {};
    \end{semilogyaxis}
    \end{tikzpicture}
    \caption{Average number of iterations}
    \label{fig:les_iterations}
\end{subfigure}
\caption{Relative error (left) and the average number of iterations required for the gradient-based optimization algorithm to converge to within a tolerance of $10^{-4}$ (right) for the turbulent airfoil problem.}
\label{fig:les_data}
\end{figure}
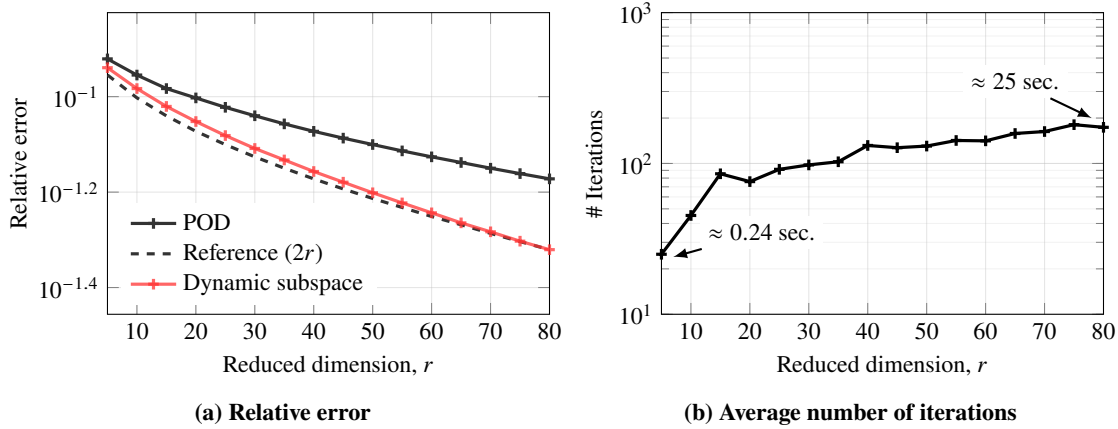

\begin{figure}[!tbp]
    \centering
    \begin{subfigure}{.33\linewidth}
        \includegraphics[width=\linewidth]{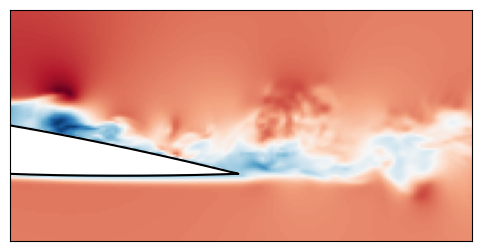}
        \caption{Exact}
    \end{subfigure}
    \begin{subfigure}{.33\linewidth}
        \includegraphics[width=\linewidth]{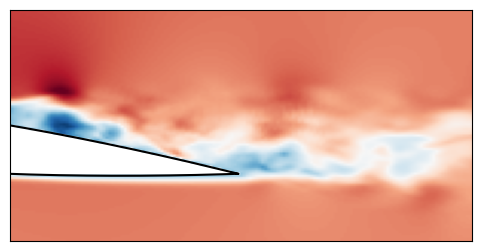}
        \caption{POD}
    \end{subfigure}
    \begin{subfigure}{.33\linewidth}
        \includegraphics[width=\linewidth]{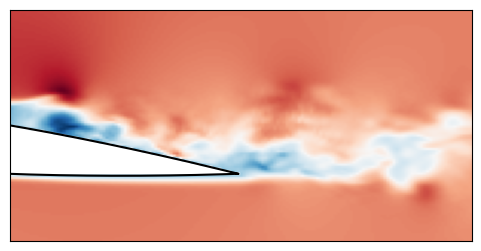}
        \caption{Dynamic subspace}
    \end{subfigure}
    \caption{Reconstruction using a dimensionality of $r=200$ at $\tau=0.25$. This particular snapshot was not included in the training data. The relative error of the entire state is 4.0\% for the POD approximation, and 2.4\% for the dynamic subspace counterpart. Reconstructions are visualized by the velocity field in the $x$-direction.}
    \label{fig:les_reconstruction}
\end{figure}

The snapshots are collected at a constant convective time increment. Each snapshot provides all velocity components extracted over the midspan slice. The full $N$ state dimension of the LES data is 666,630 with 222,210 grid points and three components of the velocity field that have been non-dimensionalized by the free stream flow velocity. For training and test data set construction, we use a total of $16{,}000$ time-resolved snapshots comprising approximately $248$ shedding cycles as explained in~\cite{towne2023database}, totaling about 43\,GB in disk storage. Rather than seeking a dynamic approximation for the entire temporal domain, we confine our analysis to the initial 10\% of the time window. This subset comprises $1{,}600$ snapshots, which are partitioned into two equally sized sets for training and for testing. The training data set consists of $K=800$ snapshots $\mathbf{Q}_\text{train} = [~\mathbf{q}(t_1)~~\mathbf{q}(t_3)~~\cdots~~\mathbf{q}(t_{1599})~] \in \mathbb{R}^{666,630 \times 800}$ to compute low-dimensional approximations. Similar to the previous example, $800$ snapshots are used to construct the test data set $\mathbf{Q}_\text{test} = [~\mathbf{q}(t_2)~~\mathbf{q}(t_4)~~\cdots~~\mathbf{q}(t_{1600})~] \in \mathbb{R}^{666,630 \times 800}$ for verification. Several representative training points, the mean flow, and the singular values of the mean-subtracted data matrix are shown in \Cref{fig:les_overview}. The number of feature modes was chosen to be $n_f = 600$.

\begin{figure}[!tbp]
    \centering
    \begin{subfigure}{.42\linewidth}
        \includegraphics[width=\linewidth]{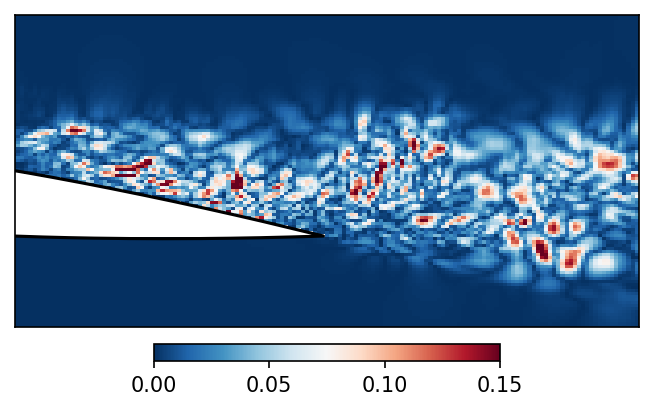}
        \caption{POD}
    \end{subfigure}
    \hspace{0.5em}
    \begin{subfigure}{.42\linewidth}
      \includegraphics[width=\linewidth]{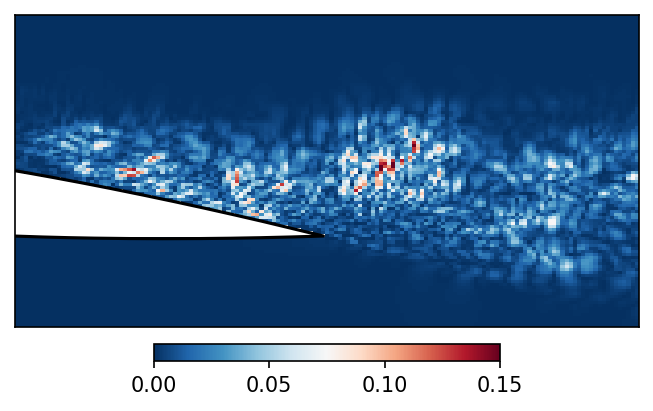}
       \caption{Dynamic subspace}
    \end{subfigure}
    \caption{Absolute error using a dimensionality of $r=200$ at $\tau=0.25$. This particular snapshot was not included in the training data.}
    \label{fig:les_absolute_error}
\end{figure}

\begin{figure}[!tbp]
    \centering
    \begin{subfigure}{.33\linewidth}
        \includegraphics[width=\linewidth]{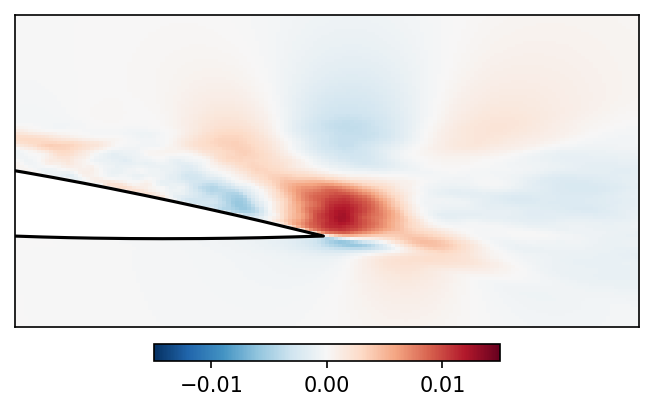}
        \caption{POD Mode \#1 (static)}
    \end{subfigure}
    \begin{subfigure}{.33\linewidth}
        \includegraphics[width=\linewidth]{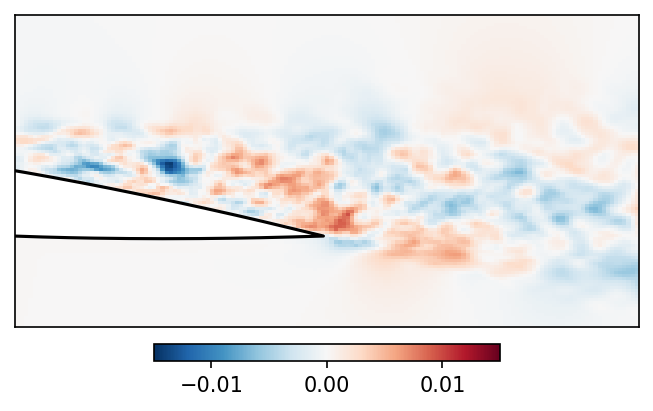}
        \caption{Dynamic mode \#1 ($\tau=0$)}
    \end{subfigure}
    \begin{subfigure}{.33\linewidth}
        \includegraphics[width=\linewidth]{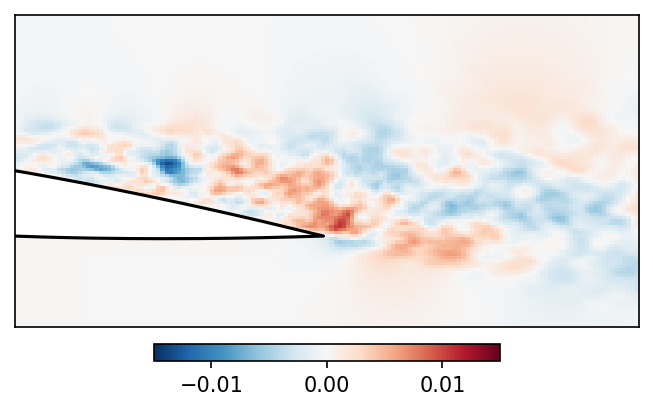}
        \caption{Dynamic Mode \#1 ($\tau=0.25$)}
    \end{subfigure}
    \begin{subfigure}{.33\linewidth}
        \includegraphics[width=\linewidth]{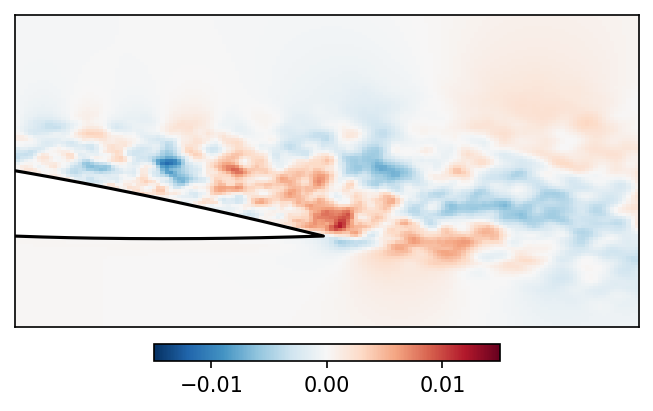}
        \caption{Dynamic Mode \#1 ($\tau=0.50$)}
    \end{subfigure}
    \begin{subfigure}{.33\linewidth}
        \includegraphics[width=\linewidth]{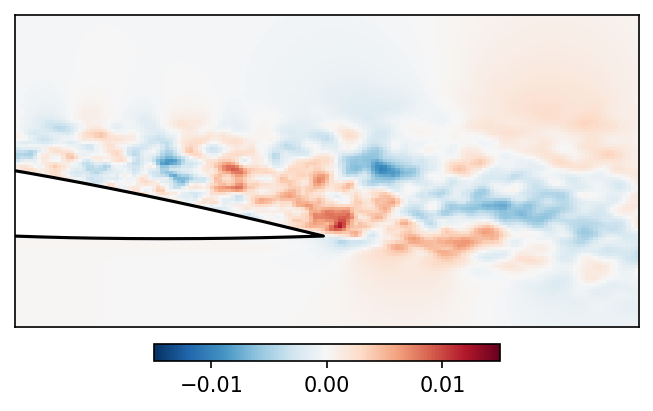}
        \caption{Dynamic Mode \#1 ($\tau=0.75$)}
    \end{subfigure}
    \begin{subfigure}{.33\linewidth}
        \includegraphics[width=\linewidth]{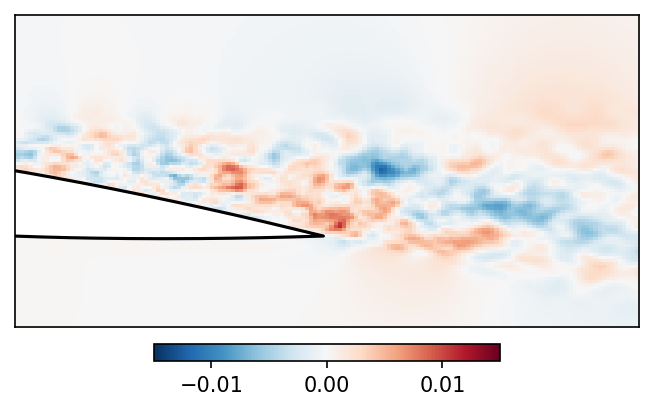}
        \caption{Dynamic Mode \#1 ($\tau=1$)}
    \end{subfigure}
    \caption{Reconstruction using a dimensionality of $r=200$ at $\tau=0.25$. This particular snapshot was not included in the training data. The relative error of the entire state is 4.0\% for the POD approximation, and 2.4\% for the dynamic subspace counterpart. POD modes are visualized by their contribution in the $x$-direction.}
    \label{fig:les_modes}
\end{figure}

\Cref{fig:les_error} illustrates the relative error for $r$-dimensional approximations of the test data. The dynamic subspace approach consistently outperforms linear approximations across all values of $r$. Similar to the 1D transport results in \Cref{subsec:one_d_transport}, this method achieves a level of precision comparable to a linear subspace with double the number of coefficients. Importantly, this accuracy gain comes with minimal computational overhead. As shown in \Cref{fig:les_iterations}, the Riemannian optimization converges to a tolerance of $10^{-4}$ in approximately 25 seconds for $r=80$. This performance was measured on an AMD Ryzen 7 PRO 8840HS CPU and represents an average across five randomized initializations of $\boldsymbol{\Theta}$. Furthermore, it should be noted that using fewer than $n_f=600$ features can further accelerate the computation with negligible loss in accuracy, provided the feature subspace remains sufficiently rich to surpass standard POD performance.

For illustration purposes, we show the reconstructions for a test snapshot at a normalized time of $t=0.25$ at a reduced dimension of $r=200$ in \Cref{fig:les_reconstruction}. We focused on a magnified region in the near-wake region as indicated in Fig~\ref{fig:les_overview}.
Qualitatively, the dynamic subspace approach captures more fine-scale detail of the turbulent flow features than the standard POD, which tends to produce smoother, more diffused results.
For this specific snapshot, the dynamic subspace approach achieves a relative error of $2.4\%$, compared to $4.0\%$ for POD. These improvements are even more pronounced in the pointwise error distributions, shown in \Cref{fig:les_absolute_error}, which show that the dynamic subspace method reduces the error magnitude in the airfoil wake and more effectively represents complex spatial variations compared to POD.

Finally, \Cref{fig:les_modes} compares the dominant POD mode and the first dynamic mode at different times. While the POD mode is inherently static and captures only the time-averaged, large-scale coherent structures, the dynamic basis vectors are more adept at resolving finer turbulent scales. More crucially, these dynamic basis vectors exhibit time continuity, allowing the basis to evolve smoothly with the flow physics and substantially enhancing the expressiveness of the reduced-order approximation. We also emphasize that, unlike other time-dependent manifold methods, our dynamic basis vectors do not require explicit orthogonalization: they remain orthonormal by design through the Riemannian geometry of the underlying manifold, representing a significant computational advantage of the proposed framework.

\section{Conclusion and Future Work}
\label{sec:conlusions_and_future_work}

This work introduces a dynamic subspace approach that effectively addresses the fundamental limitations of static linear subspaces in large-scale nonlinear systems. By parameterizing a low-dimensional basis as a geodesic on the Grassmannian manifold, the framework provides a time-continuous, adaptive basis capable of tracking non-stationary phenomena and traveling waves. This geometric formulation mitigates the Kolmogorov width barrier, which often necessitates computationally prohibitive rank inflation in traditional POD. Through the use of a feature-based reduction strategy, the optimization remains scalable and independent of the original state dimension, making it suitable for high-fidelity datasets such as the turbulent airfoil wake analyzed herein. Crucially, across both the 1D transport equation and the large-scale turbulent airfoil wake LES, the framework demonstrates the ability to recover reconstruction performance nearing that of a static basis with rank $2r$ while utilizing only $r$ degrees of freedom. This effectively doubles the representational efficiency of the reduced space and highlights the power of geodesic parameterization in capturing dominant directional shifts. Such results offer a robust path forward for modeling transport-dominated flows, providing near-optimal subspace alignment with significantly lower computational complexity.

Future research will focus on the extension to piecewise geodesics to capture complex trajectories over longer time horizons, the coupling of this adaptive basis with existing data-driven model reduction frameworks, and continued code development with the goal of open-sourcing the implementation for the broader scientific community.

\section*{Appendix: Feature Space Formulation and Vectorization}
\label{app:vectorization}

This appendix presents the projection and vectorization steps used to efficiently evaluate the optimization objective. By projecting the high-dimensional data onto the low-dimensional basis $\tilde{\mathbf{V}}$, the entire objective function \eqref{eq:reformulated} is first recast into the feature space. This isolates the invariant snapshot residual and drastically reduces the computational footprint of subsequent matrix operations. Secondly, evaluating the objective via iterative looping over each time step $\tau_j$ introduces a massive bottleneck when $K \gg 1$. To take advantage of modern parallel processing architectures, we demonstrate how this formulation is compressed into a loop-free, vectorized matrix expression.
  
\subsection*{Feature Space Formulation}

We decompose any snapshot vector $\mathbf{x}_j = \mathbf{q}(t_j) - \bar{\mathbf{q}}$ into two orthogonal components 
\begin{equation}
    \mathbf{x}_j = \mathbf{x}^{\parallel}_j + \mathbf{x}^{\perp}_j, 
\end{equation}
where $\mathbf{x}^{\parallel}_j = \tilde{\mathbf{V}}\tilde{\mathbf{V}}^\top \mathbf{x}_j$ denotes the projection of $\mathbf{x}_j$ onto the columns span of $\tilde{\mathbf{V}}$; and $\mathbf{x}^{\perp}_j = (\mathbf{I} - \tilde{\mathbf{V}}\tilde{\mathbf{V}}^\top) \mathbf{x}_j$ the residual that resides in its orthogonal complement. If we plug this decomposition into the objective function from \cref{eq:reformulated}, we obtain
\begin{equation}
    \sum_{j=1}^K \left\|
        \big(\mathbf{I} - \tilde{\mathbf{V}}\mathbf{T}\mathbf{W}(\boldsymbol{\Theta},\tau_j)\mathbf{W}(\boldsymbol{\Theta},\tau_j)^\top\mathbf{T}^\top\tilde{\mathbf{V}}^\top\big)(\mathbf{x}^{\parallel}_j  +  \mathbf{x}^{\perp}_j)
    \right\|_2^2.
\end{equation}
Because $\mathbf{x}^{\perp}_j$ is orthogonal to the columns of $\tilde{\mathbf{V}}$, it vanishes when multiplied by $\tilde{\mathbf{V}}^\top$:
\begin{equation}
    \sum_{j=1}^K \left\| \mathbf{x}^{\parallel}_j - \tilde{\mathbf{V}}\mathbf{T}\mathbf{W}(\boldsymbol{\Theta},\tau_j)\mathbf{W}(\boldsymbol{\Theta},\tau_j)^\top\mathbf{T}^\top \tilde{\mathbf{V}}^\top \mathbf{x}^{\parallel}_j + \mathbf{x}^{\perp}_j
    \right\|_2^2.
\end{equation}
Because $\big(\mathbf{x}^{\parallel}_j - \tilde{\mathbf{V}}\mathbf{T}\mathbf{W}(\boldsymbol{\Theta},\tau_j)\mathbf{W}(\boldsymbol{\Theta},\tau_j)^\top\mathbf{T}^\top \tilde{\mathbf{V}}^\top \mathbf{x}^{\parallel}_j \big)$ and $\mathbf{x}^{\perp}_j$ are orthogonal, we may split the objective into two separate norms:
\begin{equation}
    \sum_{j=1}^K \left( \left\| \mathbf{x}^{\parallel}_j - \tilde{\mathbf{V}}\mathbf{T}\mathbf{W}(\boldsymbol{\Theta},\tau_j)\mathbf{W}(\boldsymbol{\Theta},\tau_j)^\top\mathbf{T}^\top \tilde{\mathbf{V}}^\top \mathbf{x}^{\parallel}_j \right\|_2^2 + \big\| \mathbf{x}^{\perp}_j
    \big\|_2^2 \right).
\end{equation}
Plugging in the definition for $\mathbf{x}^{\parallel}_j$ and factoring out $\tilde{\mathbf{V}}$, we obtain
\begin{equation}
    \sum_{j=1}^K \left( \left\| \tilde{\mathbf{V}} \left( \tilde{\mathbf{V}}^\top \mathbf{x}_j - \mathbf{T}\mathbf{W}(\boldsymbol{\Theta},\tau_j)\mathbf{W}(\boldsymbol{\Theta},\tau_j)^\top\mathbf{T}^\top \tilde{\mathbf{V}}^\top \tilde{\mathbf{V}}\tilde{\mathbf{V}}^\top \mathbf{x}_j \right) \right\|_2^2 + \big\| \mathbf{x}^{\perp}_j
    \big\|_2^2 \right).
\end{equation}
Multiplying on the left by the orthonormal matrix $\tilde{\mathbf{V}}$ preserves the norm, which means that the above expression can fully be written in the feature space: 
\begin{equation}
    \sum_{j=1}^K \left( \left\| \tilde{\mathbf{V}}^\top \mathbf{x}_j - \mathbf{T}\mathbf{W}(\boldsymbol{\Theta},\tau_j)\mathbf{W}(\boldsymbol{\Theta},\tau_j)^\top\mathbf{T}^\top \tilde{\mathbf{V}}^\top \mathbf{x}_j \right\|_2^2 + \big\| \mathbf{x}^{\perp}_j
    \big\|_2^2 \right)
\end{equation}
We note that the term $\big\| \mathbf{x}^{\perp}_j\big\|_2^2$ is constant that does not depend on the optimal parameters $\boldsymbol{\Theta}$ and $\mathbf{T}$. It may therefore be \emph{omitted} from the optimization problem. Following the definition \cref{eq:feature_space_representation2} for representing the state data snapshots in the feature space, we then obtain
\begin{equation}
    \sum_{j=1}^K \left\|
        \mathbf{y}_j - \mathbf{T}\mathbf{W}(\boldsymbol{\Theta},\tau_j)\mathbf{W}(\boldsymbol{\Theta},\tau_j)^\top\mathbf{T}^\top \mathbf{y}_j
    \right\|_2^2.
    \label{eq:appendix1}
\end{equation}

\subsection*{Vectorized Form}

Recalling the definitions for $\mathbf{W}(\boldsymbol{\Theta},\tau_j)$ and $\mathbf{T}$ as defined in \cref{eq:omegaW} and \cref{eq:T}, respectively, the expression \cref{eq:appendix1} becomes
\begin{equation}
    \sum_{j=1}^K \left\|
        \mathbf{y}_j - \mathbf{T}\left[\begin{array}{c}
        \cos(\tau_j \boldsymbol{\Theta}) \\ \sin(\tau_j \boldsymbol{\Theta})  
    \end{array}\right]\left[~\cos(\tau_j \boldsymbol{\Theta})~~\sin(\tau_j \boldsymbol{\Theta})  
    ~\right]\left[\begin{array}{c}
        \mathbf{T}_1^\top \\ \mathbf{T}_2^\top \end{array}\right] \mathbf{y}_j
    \right\|_2^2.
    \label{eq:appendix2}
\end{equation}
After multiplying the trigonometric functions and simplifying, we have
\begin{equation}
    \sum_{j=1}^K \Bigg\|
        \mathbf{y}_j - \mathbf{T} \underbrace{\left[\begin{array}{c}
        \cos^2(\tau_j \boldsymbol{\Theta}) \mathbf{T}_1^\top \mathbf{y}_j + \cos(\tau_j \boldsymbol{\Theta}) \sin(\tau_j \boldsymbol{\Theta}) \mathbf{T}_2^\top \mathbf{y}_j \\ \sin(\tau_j \boldsymbol{\Theta}) \cos(\tau_j \boldsymbol{\Theta})\mathbf{T}_1^\top \mathbf{y}_j + \sin^2(\tau_j \boldsymbol{\Theta}) \mathbf{T}_2^\top \mathbf{y}_j
    \end{array}\right]}_{\textstyle \in \mathbb{R}^{2r}}
    \Bigg\|_2^2.
    \label{eq:appendix3}
\end{equation}
To make computing \cref{eq:appendix3} computationally efficient and avoid manual looping, we can construct the time-dependent block matrix explicitly using Hadamard products. Grouping the trigonometric functions over all time-steps into matrices where the $j$th row corresponds to the evaluation at $\tau_j$, we introduce
$\mathbf{C} = \mathbf{C}(\boldsymbol{\Theta})\in\mathbb{R}^{K\times r}$ and
$\mathbf{S} = \mathbf{S}(\boldsymbol{\Theta})\in\mathbb{R}^{K\times r}$ as
\begin{equation}
    \mathbf{C}(\boldsymbol{\Theta}) = \begin{bmatrix}
        \cos(\tau_1\theta_1) & \cos(\tau_1\theta_2) & \dots & \cos(\tau_1\theta_r) \\
        \cos(\tau_2\theta_1) & \cos(\tau_2\theta_2) & \dots & \cos(\tau_2\theta_r) \\
        \vdots &  & \ddots & \vdots \\
        \cos(\tau_K\theta_1) & \cos(\tau_K\theta_2) &\dots & \cos(\tau_K\theta_r) 
    \end{bmatrix},
    \quad
    ~\mathbf{S}(\boldsymbol{\Theta}) = \begin{bmatrix}
        \sin(\tau_1\theta_1) & \sin(\tau_1\theta_2) & \dots & \sin(\tau_1\theta_r) \\
        \sin(\tau_2\theta_1) & \sin(\tau_2\theta_2) & \dots & \sin(\tau_2\theta_r) \\
        \vdots &  & \ddots & \vdots \\
        \sin(\tau_K\theta_1) & \sin(\tau_K\theta_2) &\dots & \sin(\tau_K\theta_r) 
    \end{bmatrix}.
\end{equation}
In the Frobenius norm, the vectorized objective then takes the form
\begin{equation}
    \left\|
        \mathbf{Y} - \mathbf{T} \mathbf{R}(\mathbf{T, \boldsymbol{\Theta}},\mathbf{Y})
    \right\|_F^2
    \label{eq:appendix4}
\end{equation}
with $\mathbf{R}(\mathbf{T, \boldsymbol{\Theta}},\mathbf{Y})$ defined in accord with \cref{eq:R}.

\section*{Acknowledgments}
The authors thank Laura Balzano from the University of Michigan for several fruitful discussions on dynamic subspace estimation and Riemannian optimization.

\bibliography{bib}

@article{babaee2016minimization,
    title = {A minimization principle for the description of modes associated with finite-time instabilities},
    author = {Babaee, H and Sapsis, T P},
    journal = {Proceedings of the Royal Society A: Mathematical, Physical and Engineering Sciences},
    publisher = {The Royal Society Publishing},
    volume = {472},
    number = {2186},
    pages = {20150779},
    year = {2016},
    doi = {10.1098/rspa.2015.0779},
}

@article{barnett2022quadratic,
    title = {Quadratic approximation manifold for mitigating the Kolmogorov barrier in nonlinear projection-based model order reduction},
    author = {Barnett, Joshua and Farhat, Charbel},
    journal = {Journal of Computational Physics},
    publisher = {Elsevier},
    volume = {464},
    pages = {111348},
    year = {2022},
    doi = {10.1016/j.jcp.2022.111348},
}

@article{benner2015pmorsurvey,
    title = {A survey of projection-based model reduction methods for parametric dynamical systems},
    author = {Benner, Peter and G{\"u}{\u{g}}ercin, Serkan and Willcox, Karen},
    journal = {SIAM Review},
    publisher = {SIAM},
    volume = {57},
    number = {4},
    pages = {483--531},
    year = {2015},
    doi = {10.1137/130932715},
}

@misc{blocker2023dynamicsubspaceestimationgrassmannian,
    title = {Dynamic subspace estimation with {G}rassmannian geodesics},
    author = {Cameron J. Blocker and Haroon Raja and Jeffrey A. Fessler and Laura Balzano},
    year = {2023},
    howpublished = {arXiv:2303.14851},
    doi = {10.48550/arXiv.2303.14851},
}

@misc{deepblue,
    title = {{Turbulent airfoil wake large eddy simulation [Data set], University of Michigan - Deep Blue Data}},
    author = {Towne, A. and Yeh, C. and Patel, H. and Taira, K.},
    year = {2022},
    howpublished = {\url{http://deepblue.lib.umich.edu/data/collections/kk91fk98z}},
}

@article{edelman1998geometry,
    title = {The geometry of algorithms with orthogonality constraints},
    author = {Edelman, Alan and Arias, Tom{\'a}s A and Smith, Steven T},
    journal = {SIAM Journal on Matrix Analysis and Applications},
    publisher = {SIAM},
    volume = {20},
    number = {2},
    pages = {303--353},
    year = {1998},
    doi = {10.1137/S0895479895290954},
}

@article{absil2004riemannian,
    title = {Riemannian geometry of {G}rassmann manifolds with a view on algorithmic computation},
    author = {Absil, P-A and Mahony, Robert and Sepulchre, Rodolphe},
    journal = {Acta Applicandae Mathematica},
    publisher = {Springer},
    volume = {80},
    number = {2},
    pages = {199--220},
    year = {2004},
    doi = {10.1023/B:ACAP.0000013855.14971.91},
}

@article{geelen2023operator,
    title = {Operator inference for non-intrusive model reduction with quadratic manifolds},
    author = {Geelen, Rudy and Wright, Stephen and Willcox, Karen},
    journal = {Computer Methods in Applied Mechanics and Engineering},
    publisher = {Elsevier},
    volume = {403},
    pages = {115717},
    year = {2023},
    doi = {10.1016/j.cma.2022.115717},
}

@article{geelen2024learning,
    title = {Learning physics-based reduced-order models from data using nonlinear manifolds},
    author = {Geelen, Rudy and Balzano, Laura and Wright, Stephen and Willcox, Karen},
    journal = {Chaos: An Interdisciplinary Journal of Nonlinear Science},
    publisher = {AIP Publishing},
    volume = {34},
    number = {3},
    year = {2024},
    doi = {10.1063/5.0170105},
}

@article{yeh2019resolvent,
    title = {Resolvent-analysis-based design of airfoil separation control},
    author = {Yeh, Chi-An and Taira, Kunihiko},
    journal = {Journal of Fluid Mechanics},
    publisher = {Cambridge University Press},
    volume = {867},
    pages = {572--610},
    year = {2019},
    doi = {10.1017/jfm.2019.163},
}

@article{townsend2016pymanopt,
    title = {Pymanopt: {A} {P}ython toolbox for optimization on manifolds using automatic differentiation},
    author = {Townsend, James and Koep, Niklas and Weichwald, Sebastian},
    journal = {Journal of Machine Learning Research},
    volume = {17},
    number = {137},
    pages = {1--5},
    year = {2016},
    url = {http://jmlr.org/papers/v17/16-177.html},
}

@misc{jax2018github,
    title = {{JAX}: {C}omposable transformations of {P}ython+{N}um{P}y programs},
    author = {James Bradbury and Roy Frostig and Peter Hawkins and Matthew James Johnson and Yash Katariya and Chris Leary and Dougal Maclaurin and George Necula and Adam Paszke and Jake Vander{P}las and Skye Wanderman-{M}ilne and Qiao Zhang},
    year = {2018},
    howpublished = {GitHub:jax-ml/jax v0.3.13},
    url = {http://github.com/jax-ml/jax},
    version = {0.3.13},
}

@misc{hesthaven2026nonlinearmodelreductiontransportdominated,
    title = {Nonlinear model reduction for transport-dominated problems},
    author = {Jan S. Hesthaven and Benjamin Peherstorfer and Benjamin Unger},
    year = {2026},
    howpublished = {arXiv:2602.01397},
    doi = {10.48550/arXiv.2602.01397}
}

@article{lee2020model,
    title = {Model reduction of dynamical systems on nonlinear manifolds using deep convolutional autoencoders},
    author = {Lee, Kookjin and Carlberg, Kevin T},
    journal = {Journal of Computational Physics},
    publisher = {Elsevier},
    volume = {404},
    pages = {108973},
    year = {2020},
    doi = {10.1016/j.jcp.2019.108973},
}

@article{towne2023database,
    title = {A database for reduced-complexity modeling of fluid flows},
    author = {Towne, Aaron and Dawson, Scott TM and Br{\`e}s, Guillaume A and Lozano-Dur{\'a}n, Adri{\'a}n and Saxton-Fox, Theresa and Parthasarathy, Aadhy and Jones, Anya R and Biler, Hulya and Yeh, Chi-An and Patel, Het D and others},
    journal = {AIAA Journal},
    publisher = {American Institute of Aeronautics and Astronautics},
    volume = {61},
    number = {7},
    pages = {2867--2892},
    year = {2023},
    doi = {10.2514/1.J062203},
}

@article{sirovich1987turbulence,
    title = {{Turbulence and the dynamics of coherent structures. I. Coherent structures}},
    author = {Sirovich, Lawrence},
    journal = {Quarterly of Applied Mathematics},
    volume = {45},
    number = {3},
    pages = {561--571},
    year = {1987},
    doi = {10.1090/qam/910462},
}

@article{berkooz1993proper,
    title = {The proper orthogonal decomposition in the analysis of turbulent flows},
    author = {Berkooz, Gal and Holmes, Philip and Lumley, John L},
    journal = {Annual Review of Fluid Mechanics},
    publisher = {Annual Reviews 4139 El Camino Way, PO Box 10139, Palo Alto, CA 94303-0139, USA},
    volume = {25},
    number = {1},
    pages = {539--575},
    year = {1993},
    doi = {10.1146/annurev.fl.25.010193.002543},
}

@article{hesthaven2022reduced,
    title = {Reduced basis methods for time-dependent problems},
    author = {Hesthaven, Jan S and Pagliantini, Cecilia and Rozza, Gianluigi},
    journal = {Acta Numerica},
    publisher = {Cambridge University Press},
    volume = {31},
    pages = {265--345},
    year = {2022},
    doi = {10.1017/S0962492922000058},
}

@article{iollo2014advection,
    title = {Advection modes by optimal mass transfer},
    author = {Iollo, Angelo and Lombardi, Damiano},
    journal = {Physical Review E},
    publisher = {APS},
    volume = {89},
    number = {2},
    pages = {022923},
    year = {2014},
    doi = {10.1103/PhysRevE.89.022923},
}

@article{gerbeau2014approximated,
    title = {Approximated {L}ax pairs for the reduced order integration of nonlinear evolution equations},
    author = {Gerbeau, Jean-Fr{\'e}d{\'e}ric and Lombardi, Damiano},
    journal = {Journal of Computational Physics},
    publisher = {Elsevier},
    volume = {265},
    pages = {246--269},
    year = {2014},
    doi = {10.1016/j.jcp.2014.01.047},
}

@article{rapun2010reduced,
    title = {Reduced order models based on local {POD} plus {G}alerkin projection},
    author = {Rap{\'u}n, Mar{\'\i}a-Luisa and Vega, Jos{\'e} M},
    journal = {Journal of Computational Physics},
    publisher = {Elsevier},
    volume = {229},
    number = {8},
    pages = {3046--3063},
    year = {2010},
    doi = {10.1016/j.jcp.2009.12.029},
}

@article{rim2023manifold,
    title = {Manifold approximations via transported subspaces: {M}odel reduction for transport-dominated problems},
    author = {Rim, Donsub and Peherstorfer, Benjamin and Mandli, Kyle T},
    journal = {SIAM Journal on Scientific Computing},
    publisher = {SIAM},
    volume = {45},
    number = {1},
    pages = {A170--A199},
    year = {2023},
    doi = {10.1137/20M1316998},
}

@article{black2020projection,
    title = {Projection-based model reduction with dynamically transformed modes},
    author = {Black, Felix and Schulze, Philipp and Unger, Benjamin},
    journal = {ESAIM: Mathematical Modelling and Numerical Analysis},
    publisher = {EDP Sciences},
    volume = {54},
    number = {6},
    pages = {2011--2043},
    year = {2020},
    doi = {10.1051/m2an/2020046},
}

@article{peherstorfer2022breaking,
    title = {Breaking the {K}olmogorov barrier with nonlinear model reduction},
    author = {Peherstorfer, Benjamin},
    journal = {Notices of the American Mathematical Society},
    volume = {69},
    number = {5},
    pages = {725--733},
    year = {2022},
    doi = {10.1090/noti2475},
}

@article{peherstorfer2015dynamic,
    title = {Dynamic data-driven reduced-order models},
    author = {Peherstorfer, Benjamin and Willcox, Karen},
    journal = {Computer Methods in Applied Mechanics and Engineering},
    publisher = {Elsevier},
    volume = {291},
    pages = {21--41},
    year = {2015},
    doi = {10.1016/j.cma.2015.03.018},
}

@article{peherstorfer2015online,
    title = {Online adaptive model reduction for nonlinear systems via low-rank updates},
    author = {Peherstorfer, Benjamin and Willcox, Karen},
    journal = {SIAM Journal on Scientific Computing},
    publisher = {SIAM},
    volume = {37},
    number = {4},
    pages = {A2123--A2150},
    year = {2015},
    doi = {10.1137/140989169},
}

@article{billaud2017dynamical,
    title = {Dynamical model reduction method for solving parameter-dependent dynamical systems},
    author = {Billaud-Friess, Marie and Nouy, Anthony},
    journal = {SIAM Journal on Scientific Computing},
    publisher = {SIAM},
    volume = {39},
    number = {4},
    pages = {A1766--A1792},
    year = {2017},
    doi = {10.1137/16M1071493},
}

@article{koch2007dynamical,
    title = {Dynamical low-rank approximation},
    author = {Koch, Othmar and Lubich, Christian},
    journal = {SIAM Journal on Matrix Analysis and Applications},
    publisher = {SIAM},
    volume = {29},
    number = {2},
    pages = {434--454},
    year = {2007},
    doi = {10.1137/050639703},
}

@article{sapsis2009dynamically,
    title = {Dynamically orthogonal field equations for continuous stochastic dynamical systems},
    author = {Sapsis, Themistoklis P and Lermusiaux, Pierre FJ},
    journal = {Physica D: Nonlinear Phenomena},
    publisher = {Elsevier},
    volume = {238},
    number = {23-24},
    pages = {2347--2360},
    year = {2009},
    doi = {10.1016/j.physd.2009.09.017},
}

@article{ramezanian2021fly,
    title = {On-the-fly reduced order modeling of passive and reactive species via time-dependent manifolds},
    author = {Ramezanian, Donya and Nouri, Arash G and Babaee, Hessam},
    journal = {Computer Methods in Applied Mechanics and Engineering},
    publisher = {Elsevier},
    volume = {382},
    pages = {113882},
    year = {2021},
    doi = {10.1016/j.cma.2021.113882},
}

@article{patil2020real,
    title = {Real-time reduced-order modeling of stochastic partial differential equations via time-dependent subspaces},
    author = {Patil, Prerna and Babaee, Hessam},
    journal = {Journal of Computational Physics},
    publisher = {Elsevier},
    volume = {415},
    pages = {109511},
    year = {2020},
    doi = {10.1016/j.jcp.2020.109511},
}

@article{yang1995projection,
    title = {Projection approximation subspace tracking},
    author = {Yang, Bin},
    journal = {IEEE Transactions on Signal Processing},
    publisher = {IEEE},
    volume = {43},
    number = {1},
    pages = {95--107},
    year = {1995},
    doi = {10.1109/78.365290},
}

@article{miller1981moving,
    title = {Moving finite elements. I},
    author = {Miller, Keith and Miller, Robert N},
    journal = {SIAM Journal on Numerical Analysis},
    publisher = {SIAM},
    volume = {18},
    number = {6},
    pages = {1019--1032},
    year = {1981},
    doi = {10.1137/0718070},
}

@article{gugercin2004balancedtruncation,
    title = {A survey of model reduction by balanced truncation and some new results},
    author = {G{\"u}{\u{g}}ercin, Serkan and Antoulas, Athanasios C},
    journal = {International Journal of Control},
    volume = {77},
    number = {8},
    pages = {748--766},
    year = {2004},
    publisher = {Taylor \& Francis},
    doi = {10.1080/00207170410001713448},
}

@article{gosea2022datadrivenbt,
    title = {Data-driven balancing of linear dynamical systems},
    author = {Gosea, Ion Victor and G{\"u}{\u{g}}ercin, Serkan and Beattie, Christopher},
    journal = {SIAM Journal on Scientific Computing},
    volume = {44},
    number = {1},
    pages = {A554--A582},
    year = {2022},
    publisher = {SIAM},
    doi = {10.1137/21M1411081},
}

@misc{schwerdtner2024greedyquad,
    title = {Greedy construction of quadratic manifolds for nonlinear dimensionality reduction and nonlinear model reduction}, 
    author = {Paul Schwerdtner and Benjamin Peherstorfer},
    year = {2024},
    howpublished = {arXiv:2403.06732},
    doi = {10.48550/arXiv.2403.06732},
}

@inproceedings{scholkopf1997kernel,
    title = {Kernel principal component analysis},
    author = {Sch{\"o}lkopf, Bernhard and Smola, Alexander and M{\"u}ller, Klaus-Robert},
    booktitle = {International Conference on Artificial Neural Networks},
    pages = {583--588},
    year = {1997},
    doi = {10.1007/BFb0020217},
}

@article{diaz2025kernelmanifold,
    author = {Diaz, Alejandro N. and Needels, Jacob T. and Tezaur, Irina K. and Blonigan, Patrick J.},
    title = {Kernel manifolds: {N}onlinear-augmentation dimensionality reduction using reproducing kernel {H}ilbert spaces},
    journal = {International Journal for Numerical Methods in Engineering},
    volume = {126},
    number = {24},
    pages = {e70230},
    year = {2025},
    doi = {https://doi.org/10.1002/nme.70230},
}

@article{fries2022lasdi,
    title = {{LaSDI}: {P}arametric latent space dynamics identification},
    author = {William D. Fries and Xiaolong He and Youngsoo Choi},
    journal = {Computer Methods in Applied Mechanics and Engineering},
    volume = {399},
    pages = {115436},
    year = {2022},
    doi = {10.1016/j.cma.2022.115436},
}

\end{document}